\documentclass[12pt,reqno]{amsart}
\usepackage{epsfig, amsmath, amssymb,latexsym}
\usepackage{hyperref}
 \makeatletter
\def\LaTeX{\leavevmode L\raise.42ex
\hbox{\kern-.3em\size{\sf@size}{0pt}\selectfont A}\kern-.15em\TeX}
\def\@currentlabel{2.1}\label{e:dispaa}
\def\@currentlabel{2.21}\label{e:dispau}
\def\@currentlabel{2.22}\label{e:dispav}
\def\@currentlabel{2.23}\label{e:dispaw}
\def\@currentlabel{2.24}\label{e:dispax}
\def\theequation{\thesection.\@arabic\c@equation}
\makeatother
\newcommand{\R}{{\mathbb{R}}}
\newcommand{\C}{{\mathcal{C}}}
\newcommand{\V}{{\mathcal{V}}}

\renewcommand{\theequation}{\thesection.\arabic{equation}}
\newtheorem{teo}{Theorem}[section]
\newtheorem{definition}[teo]{Definition}
\newtheorem{prop}[teo]{Proposition}
\newtheorem{lema}[teo]{Lemma}

\newtheorem{remark}[teo]{Remark}

\begin{document}

\title[The square root of the Laplacian]{Positive solutions of nonlinear problems
involving the square root of the Laplacian}


\author{Xavier Cabr\'{e}}

\address{X.C., ICREA and Departament de Matem\`{a}tica
Aplicada~I, Universitat Polit\`{e}cnica de Catalunya, Diagonal 647,
08028 Barcelona, Spain}

\email{xavier.cabre@upc.edu}

\author{Jinggang Tan}

\address{J.T., Departamento de  Matem\'atica,
Universidad T\'{e}cnica Federico Santa Mar\'{i}a,  Avda. Espa\~na 1680,
Valpara\'{\i}so, Chile} 

\email{jinggang.tan@usm.cl}

 \subjclass{}
\begin{abstract}
We consider nonlinear elliptic problems involving a nonlocal operator:
the square root of the Laplacian in a bounded domain
with zero Dirichlet boundary conditions. For positive solutions
to problems with power nonlinearities, we
establish existence and regularity results, as well as
a priori estimates of Gidas-Spruck type. In addition, among
other results, we prove a symmetry theorem of Gidas-Ni-Nirenberg type.
\end{abstract}
\date{}\maketitle

\setcounter{equation}{0}
\section{Introduction}

This paper is concerned with the study of positive solutions to nonlinear
problems involving a nonlocal positive operator: the square root of
the Laplacian in a bounded domain with zero Dirichlet boundary
conditions. We look for solutions to the nonlinear problem 
\begin{equation}\label{eqn-frac}
 \left\{
\begin{array}{ll} A_{1/2}u=f(u) &\mbox{in} \; \Omega,\\
u=0  &\mbox{on}\;\partial\Omega,\\
u>0  &\mbox{in} \; \Omega,
\end{array}
\right.
\end{equation}
where $\Omega$ is a smooth bounded domain of $\R^{n}$ and $A_{1/2}$
stands for the square root of the Laplacian operator $-\Delta$ in
$\Omega$ with zero Dirichlet boundary values on $\partial \Omega$.

To define $A_ {1/2}$, let $\{\lambda_{k},\varphi_{k}\}_{k=1}^{\infty}$ be
 the eigenvalues and corresponding eigenfunctions of the Laplacian operator $-\Delta$ in
$\Omega$ with zero Dirichlet boundary values on $\partial \Omega$,
\begin{equation*}
\left\{
\begin{array}{ll}
-\Delta \varphi_{k}=\lambda_{k}\varphi_{k}& \mbox{in} \;\Omega,\\
\varphi_{k}=0&\mbox{on}\; \Omega,
\end{array}
\right.
\end{equation*}
normalized by $\|\varphi_{k}\|_{L^{2}(\Omega)}=1$.
The square root of the Dirichlet Laplacian, that we denote by
$A_{1/2}$, is given by
\begin{equation}
\label{defAfirst}
u=\sum_{k=1}^{\infty}c_{k} \varphi_{k}\longmapsto
A_{1/2}u=\sum_{k=1}^{\infty}c_{k}\lambda_{k}^{1/2}\varphi_{k},
\end{equation}
which clearly maps
$H_{0}^{1}(\Omega)=\{u=\sum_{k=1}^{\infty}c_{k} \varphi_{k}\mid\sum_{k=1}^{\infty}\lambda_{k}c_{k}^{2}<\infty
\}$ into $L^{2}(\Omega)$.

The fractions of the Laplacian, such as the previous square root
$A_{1/2}$, are the infinitesimal generators of  L{\'e}vy stable diffusion processes and
appear in anomalous diffusions in plasmas, flames propagation and
chemical reactions in liquids, population dynamics, geophysical fluid
dynamics, and American options in finance.

Essential to the results in
this paper is to realize the nonlocal operator $A_{1/2}$ in  the
following local manner. Given a function
$u$ defined in $\Omega$, we consider its harmonic extension $v$ in the cylinder
$\C:=\Omega\times(0,\infty)$, with $v$
vanishing on the lateral
boundary $\partial_{L}\C:=\partial\Omega\times[0,\infty)$.  Then,
$A_{1/2}$ is given by the Dirichlet to Neumann map on $\Omega$,
$u\mapsto \frac{\partial v}{\partial \nu}\!\!\mid_{\Omega\times\{0\}}$,
 of such
harmonic extension in the cylinder.
 In this way, we transform problem (\ref{eqn-frac}) to a local problem in one more dimension.
By studying this problem with classical local techniques, we establish existence
  of positive solutions for problems with subcritical
power nonlinearities, regularity and an $L^{\infty}$-estimate of
Brezis-Kato type for weak solutions, a priori estimates of
Gidas-Spruck type, and a nonlinear Liouville
 type result for the square root of  the Laplacian in the half-space.
 We also obtain a symmetry theorem of Gidas-Ni-Nirenberg type.

The analogue problem to (\ref{eqn-frac}) for the Laplacian has been
investigated widely in the last decades. This is the
problem
\begin{equation}\label{eqn-sub-f1}
\begin{cases}
-\Delta u =f(u) &\text{in}\; \Omega,\\
u=0& \text{on} \;\partial \Omega, \\
u>0 & \text{in} \; \Omega;
\end{cases}
\end{equation}
see~\cite{St} and references therein. Considering the minimization problem
$\text{min}\{\|u\|_{H^{1}_{0}(\Omega)}\mid
\|u\|_{L^{p+1}(\Omega)}=1\}$, one obtains a positive solution of (\ref{eqn-sub-f1}) in the
case $f(u)=u^{p}$, $1<p<\frac{n+2}{n-2}$, since the Sobolev
embedding is compact. Ambrosetti and Rabinowitz \cite{AR} introduced
the mountain pass theorem to study problem (\ref{eqn-sub-f1}) for
more general subcritical nonlinearities. Instead, for
$f(u)=u^{\frac{n+2}{n-2}}$, Pohozaev identity leads to nonexistence
to (\ref{eqn-sub-f1}) if $\Omega$ is star-shaped.
In contrast, Brezis and Nirenberg \cite{BN} showed that the nonexistence
of solution may be reverted by adding  a small linear perturbation to the critical
power nonlinearity.

For the square root  $A_{1/2}$ of the Laplacian, we derive the
following result on existence of positive solutions to problem
(\ref{eqn-frac}).

\begin{teo}\label{teo-sub1} Let $n\ge 1$ be an integer and
 $2^{\sharp}=\frac{2n}{n-1}$ when $n\ge 2$.
Suppose that $\Omega$ is a smooth bounded domain in $\R^{n}$ and
$f(u)=u^{p}$.  Assume that $1<p<2^{\sharp}-1=\frac{n+1}{n-1}$ if
$n\ge 2$, or $1< p<\infty$ if $n=1$.

Then,  problem
{\rm(\ref{eqn-frac})} admits at least one solution. This solution
$($as well as every weak solution$)$ belongs to
$C^{2,\alpha}(\overline{\Omega})$ for some $0<\alpha<1$.
\end{teo}

As mentioned before, we realize problem (\ref{eqn-frac}) through
a local problem in one more dimension by a
Dirichlet to Neumann map. This provides a variational
structure to the problem, and we study its corresponding minimization problem.
Here the Sobolev trace embedding comes into play, and its critical exponent
$2^{\sharp}=\frac{2n}{n-1}$, $n\ge 2$, is the power appearing in Theorem \ref{teo-sub1}.
We call $p$ critical (respectively, subcritical or supercritical) when
$p=2^{\sharp}-1=\frac{n+1}{n-1}$ (respectively, $p<2^{\sharp}-1$ or  $p>2^{\sharp}-1$). 
In the subcritical case of Theorem
\ref{teo-sub1}, the compactness of the Sobolev trace embedding  in bounded  domains leads to
the existence of solution. Its regularity will be consequence of further results
presented later in this introduction.

\begin{remark}{\rm In \cite{Tan} the second author J. Tan establishes  the non-existence
of classical solutions to (\ref{eqn-frac}) with
$f(u)=u^{p}$ in star-shaped domains for the critical and  supercritical
 cases. In addition, an existence result of
  Brezis-Nirenberg type \cite{BN} for $f(u)=u^{p}+\mu u$, $\mu>0$, is also established.}
\end{remark}

 Gidas and Spruck \cite{GS}
established  a priori estimates for positive solutions of problem
(\ref{eqn-sub-f1}) when $f(u)=u^{p}$ and $p<\frac{n+2}{n-2}$. Its
proof involves the method of blow-up combined with two
important ingredients: nonlinear Liouville type results in all space
and in a half-space. The proofs of such Liouville theorems are based on
 the Kelvin transform and the moving planes method or the moving spheres method.
 Here we establish an analogue:  the following  a priori estimates of Gidas-Spruck
type for solutions of problem (\ref{eqn-frac}).

\begin{teo} \label{teo-ap} Let $n\ge 2$ and  $2^{\sharp}=\frac{2n}{n-1}$.
Assume that $\Omega\subset \R^{n}$ is a smooth bounded  domain and
$f(u)=u^{p}$, $1< p<2^{\sharp}-1=\frac{n+1}{n-1}$.

Then,  there
exists a constant $C(p,\Omega)$, which depends only on $p$ and~$\Omega$, 
such that every weak solution of {\rm(\ref{eqn-frac})} satisfies
\[
\|u\|_{L^{\infty}(\Omega)}\le C(p,\Omega).
\]
\end{teo}

To prove this result, we combine the blow-up method and two useful
ingredients: a nonlinear Liouville theorem for the square root of the
Laplacian in all of  $\R^{n}$, and a similar one in the half-space
$\R^{n}_{+}$ with zero Dirichlet boundary values
on $\partial \R^{n}_{+}$. The first one in the whole space was proved by
Ou \cite{Ou} using the moving planes method and  by Y.Y.
Li, M. Zhu and L. Zhang \cite{LiZhu95}, \cite{LiZ} using the moving spheres method.
Its statement is the following.

\begin{teo} {\rm(}{\bf\cite{LiZhu95}}, {\bf\cite{Ou}}, {\bf\cite{LiZ}}{\rm)} \label{teo-LZO-half}
 For $n\ge 2$ and  $1< p<2^{\sharp}-1=\frac{n+1}{n-1}$, there exists no weak solution of problem
\begin{equation}\label{eqn-half-Liou}
\begin{cases}
(-\Delta)^{1/2}u=u^{p} &\text{in} \; \R^{n},\\
u >0  &\text{in}\; \R^{n}.
\end{cases}
\end{equation}
\end{teo}
 As we will see later, here $(-\Delta)^{1/2}$ is the usual half-Laplacian in all of $\R^{n}$, and
  problem (\ref{eqn-half-Liou}) is equivalent to problem $\Delta v=0$ and $v>0$ in $\R^{n+1}_{+}$,
   $\partial_{\nu}v=v^{p}$ on $\partial \R^{n+1}_{+}$. The corresponding Liouville theorem
for the square root of the Laplacian 
in $\R_{+}^{n}=\{x\in \R^{n} \mid x_{n}>0\}$ was not available and we establish it in this
paper for bounded solutions.
\begin{teo} \label{thm-l-d2} Let $n\ge 2$,
$2^{\sharp}=\frac{2n}{n-1}$, and $1<p\le
2^{\sharp}-1=\frac{n+1}{n-1}$. Then, there exists no bounded solution $u$
of
\begin{equation}
 \left\{
\begin{array}{ll} A_{1/2}u=u^{p}  &\mbox{in} \; \R^{n}_{+},\\
u=0  &\mbox{on}\;\partial\R^{n}_{+},\\
u>0  &\mbox{in} \; \R^{n}_{+},
\end{array}
\right.
\end{equation}
where $A_{1/2}$ is the square root of the Laplacian in $\R_{+}^{n}=\{ x_{n}>0\}$
with zero Dirichlet boundary conditions on $\partial \R^{n}_{+}$.

In an equivalent way, let
\[
\R^{n+1}_{++}=\{z=(x_{1},x_{2},\cdots,x_{n},y)\mid x_{n}>0, y > 0\}.
\]
 If $n\ge 2$ and  $1< p\le 2^{\sharp}-1=\frac{n+1}{n-1}$, then there exists no
  bounded solution $v\in C^{2}(\R^{n+1}_{++})\cap
C(\overline{\R^{n+1}_{++}})$ of
\begin{equation}
\begin{cases}
\Delta  v=0 &  \text{in}\; \R^{n+1}_{++}, \\
v=0  &\text{on}\; \{x_{n}=0, y>0\},\\
\frac{\partial v}{\partial \nu}=v^{p} & \text{on}\; \{x_{n}>0, y=0\}, \\
 v>0 &  \text{in}\;  \R^{n+1}_{++},
\end{cases}
\end{equation}
 where $\nu$ is the unit outer normal to $\R^{n+1}_{++}$ at $\{x_{n}>0, y=0\}$.
\end{teo}
The proof of this result combines the Kelvin transform, the moving
planes method, and  a Hamiltonian identity for the half-Laplacian found by
Cabr{\'e} and Sol\`{a}-Morales \cite{CS05}. The result of Theorem
\ref{thm-l-d2} is still open without the assumption of boundedness
of the solution.

Gidas, Ni, and Nirenberg  \cite{GNN} established symmetry properties
for solutions to problem (\ref{eqn-sub-f1}) when $f$ is Lipschitz
continuous and $\Omega$ has certain symmetries. The proof of these
symmetry results uses the maximum principle and the moving planes
method.  The moving planes method was introduced by Alexandroff to study
a geometric problem, while in the framework of problem (\ref{eqn-sub-f1}) was first used by Serrin.
In the improved version of Berestycki and Nirenberg \cite{BN91}, it
replaces the use of Hopf's lemma by a maximum principle in domains
of small measure.

Here we proceed in a similar manner and obtain the following symmetry result of
Gidas-Ni-Nirenberg type for (\ref{eqn-frac}).

\begin{teo}\label{teo-sys} Assume that $\Omega$ is a bounded smooth domain of  $\R^{n}$ which
is convex in the $x_{1}$ direction and symmetric with respect to the
hyperplane $\{x_{1}=0\}$. Let   $f$ be Lipschitz continuous and
$u$ be a $C^{2}(\overline{\Omega})$ solution of {\rm(\ref{eqn-frac})}.

Then, $u$ is symmetric with respect to $x_{1}$, i.e.,
$u(-x_{1},x')=u(x_{1},x')$ for all $(x_{1},x')\in\Omega$. In
addition, $\frac{\partial u}{\partial x_{1}}<0$ for $x_{1}>0$.

In particular, if $\Omega=B_{R}(0)$ is a ball, then $u$ is radially
symmetric, $u=u(|x|)=u(r)$ for $r=|x|$, and it is decreasing,
i.e., $u_{r}<0$ for $0<r<R$.
\end{teo}
We prove this symmetry result by using the moving planes method combined with
 the following maximum principle for the square root $A_{1/2}$ of the Laplacian in
domains of small measure (see
Proposition \ref{lem-max3} for a more general statement in nonsmooth domains).
\begin{prop}\label{prop-smalld0}
Assume that  $u\in
C^{2}(\overline{\Omega})$ satisfies
\begin{equation*}
\left\{
\begin{array}{ll}
A_{1/2}u +c(x)u \ge 0 &\mbox{in}\; \Omega,\\
u=0 &\mbox{on}\; \partial \Omega,
\end{array}
\right.
\end{equation*}
where $\Omega$ is a smooth bounded domain in $\R^{n}$ and  $c\in
L^{\infty}(\Omega)$. Then, there exists $\delta>0$ depending only on
$n$ and $\|c^{-}\|_{L^{\infty}(\Omega)}$, such that if $|\Omega\cap\{u<0\}|\le
\delta$ then $u\ge 0$ in $\Omega$.
\end{prop}

The above maximum principle in ``small'' domains replaces the use of
Hopf's lemma to prove symmetry results for $A_{1/2}$ in Lipschitz
domains. We point out that Chipot, Chleb\'{i}k, Fila, and Shafrir
\cite{CCFS} studied a related problem:
\begin{equation}\label{eqn-ball}
\left\{
\begin{array}{ll}
 -\Delta v=g(v) &  \mbox{in}\; B_{R}^{+}=\{z\in \R^{n+1}\mid |z|\le R, z_{n+1}>0\},\\
 v=0 & \mbox{on}\; \partial B_{R}^{+}\cap \{z_{n+1}>0\},\\
\frac{\partial v}{\partial \nu}=f(v) & \mbox{on} \; \partial B_{R}^{+}\cap \{z_{n+1}=0\},\\
v>0 &  \mbox{in}\;B_{R}^{+},
\end{array}
\right.
\end{equation}
where $f, g\in C^{1}(\R)$ and  $\nu$ is the unit outer normal. They
proved  existence, non-existence, and axial symmetry results for
solutions of (\ref{eqn-ball}). Following one of their proofs, we establish
Hopf's lemma  for
 $A_{1/2}$, Lemma \ref{lem-hopf} below. Finally, let us mention that singular solutions and extremal
solutions of similar problems to (\ref{eqn-ball}) have been considered by
Davila, Dupaigne, and Montenegro \cite{JD1}, \cite{DD}.

As we mentioned, 
crucial to our results is  that $A_{1/2}$ is a nonlocal operator in $\Omega$ but which can be
realized through a local problem in $\Omega\times(0,\infty)$. To explain this, let us start with
 the square root of the Laplacian (or half-Laplacian) in $\R^{n}$.  Let $u$ be a
bounded continuous function in all of  $\R^{n}$. There is a unique
harmonic extension $v$ of $u$ in the half-space
$\R^{n+1}_{+}=\R^{n}\times(0,\infty)$. That is,
\begin{equation*}
\begin{cases}
\Delta v=0& \text{in}\; \R_{+}^{n+1}=\{(x,y)\in \R^{n}\times(0,\infty)\},\\
v=u& \text{on}\; \R^{n}=\partial \R_{+}^{n+1}.
\end{cases}
\end{equation*}
Consider the operator $T: u\mapsto -\partial_{y}v(\cdot,0)$. Since
$\partial_{y}v$ is  still a harmonic function, if we apply the
operator $T$ twice, we obtain
\[
(T\circ T)u=\partial_{yy}v\mid_{y=0}=-\Delta_{x} v\mid_{y=0}=-\Delta
u\;\;\;\text{in}\; \R^{n}.
\]
Thus, we see that the operator $T$ mapping the Dirichlet data
$u$ to the Neumann data $-\partial_{y}v(\cdot,0)$ is actually a square root of the Laplacian.
Indeed it coincides with
the usual half-Laplacian, see \cite{La}.

Here we introduce a new analogue extension problem in a cylinder ${\mathcal C}:=\Omega\times (0,\infty)$
in one more dimension  to realize (\ref{eqn-frac}) by a local problem in $\C$.
More precisely, we look for a function $v$ with $v(\cdot,0)=u$ in $\R^{n}$
satisfying the following mixed boundary value problem in a half-cylinder:
\begin{equation}\label{eqn-frac2}
\left\{
\begin{array}{ll}
 \Delta v=0 &  \mbox{in}\; {\mathcal C}=\Omega\times (0,\infty),\\
 v=0 & \mbox{on}\; \partial_{L}{\mathcal C}:=\partial\Omega\times [0,\infty),\\
\frac{\partial v}{\partial \nu}=f(v) & \mbox{on} \; \Omega\times \{0\},\\
v>0 &  \mbox{in}\; {\mathcal C},
\end{array}
\right.
\end{equation}
where $\nu$ is the unit outer normal to  $\C$ at $\Omega\times
\{0\}$. If $v$ satisfies (\ref{eqn-frac2}), then the trace $u$ on
$\Omega\times\{0\}$ of $v$  is a solution of
problem (\ref{eqn-frac}). Indeed,
  since $\partial_{y}v$ is harmonic and also vanishes on the lateral boundary $\partial \Omega\times[0,\infty)$,
   we see as before that the Dirichlet to Neumann map $u\mapsto -\partial_{y}v(\cdot,0)$
 is the  unique positive square root $A_{1/2}$ of the Dirichlet Laplacian in $\Omega$.

The generators of  L{\'e}vy symmetric stable diffusion processes are the
fractional powers of the Laplacian  $(-\Delta)^{s}$ in all of $\R^{n}$, $0<s<1$.
Fractional Laplacians attract nowadays much interest in physics, biology, finance,
 as well as in mathematical nonlinear analysis (see \cite{dA04}). One of the few nonlinear 
results which is not recent is due to
Sugitani \cite{Sugitani75}, who proved blow up results for solutions of heat
equations $\partial_{t}u+(-\Delta)^{s}u=f(u)$ in $\R^{n}$, for all $0<s<1$. It is important to note that the
fundamental solution of the fractional heat equation has power decay (or heavy) tails,
in contrast with the exponential
decay in case of the classical heat equation. L{\'e}vy  processes have also been applied
to model American options \cite{dA04}. As recent nonlinear works for fractional diffusions, 
let us mention the following. Caffarelli and
Silvestre \cite{CaSi} have given a new local realization of the fractional
Laplacian $(-\Delta)^{s}$, for all $0<s<1$, through the
Dirichlet-Neumann map of an appropriate degenerate elliptic operator. The
regularity of the obstacle problem for the fractional powers of the
Laplacian operator was proved by Silvestre \cite{Sil}. The optimal
regularity for such Signorini problem was improved in \cite{CaSaSi}.
Moreover, the operator $(-\Delta)^{s}$ plays an important role in
the study of the quasi-geostrophic equations in geophysical fluid
dynamics; see the important recent paper \cite{Caf-Vasseur} by Caffarelli and Vasseur.
Cabr{\'e}
and Sol\`{a}-Morales \cite{CS05} studied layer solutions (solutions
which are monotone with respect to one variable) of
$(-\Delta)^{1/2}u=f(u)\; \text{in} \; \R^{n}$, where $f$ is of
balanced bistable type.

To prove Theorem \ref{teo-sub1}, in view of (\ref{eqn-frac2}) being a local realization of (\ref{eqn-frac}),
 we consider the Sobolev space
\begin{equation*}
H_{0,L}^{1}(\C)=\{v\in {H}^{1}({\mathcal{C}})\mid v=0
\;\mbox{a.e.}\;\;\mbox{on}\;
\partial_{L}{\mathcal{C}}=\partial\Omega\times[0,\infty)\,\},
\end{equation*}
equipped with the norm
$\|v\|=\left(\int_{\C}|\nabla v|^{2}\,dxdy\right)^{1/2}$.
Since problem (\ref{eqn-frac2}) has variational structure,
 we consider  its corresponding minimization problem
\[
I_{0}=\inf\left\{ \int_{{\mathcal{C}}}|\nabla v(x,y)|^{2}dxdy \mid
v\in H_{0,L}^{1}(\C), \int_{\Omega}|v(x,0)|^{p+1}dx=1\, \right\}.
\]
We will prove that, for subcritical powers, there is a minimizer for this problem. Its trace on $\Omega\times\{0\}$ will
provide with a weak solution of (\ref{eqn-frac}).

Thus, it is important to characterize  the space  $\V_{0}(\Omega)$ of all traces on
$\Omega\times\{0\}$ of functions in $H_{0,L}^{1}(\C)$. This is stated in the
following result ---which corresponds to Proposition
\ref{prop-V-def} of next section.

\begin{prop}\label{prop-secone-V0}
Let $\V_{0}(\Omega)$ be the space of all traces on
$\Omega\times\{0\}$ of functions in $H_{0,L}^{1}(\C)$. Then, we
have
\begin{equation*}
\begin{split}
\V_{0}(\Omega)& :=\left\{u =\mbox{{\rm tr}}_{\Omega}v\mid v\in
H_{0,L}^{1}(\C)\right\}
\\
&=\left\{u\in H^{1/2}(\Omega)\mid\int_{\Omega}\frac{u^{2}(x)}{d(x)}\,dx<+\infty\right\}\\
&= \left\{ u \in L^{2}(\Omega) \mid
u=\sum_{k=1}^{\infty}b_{k}\varphi_{k} \;\text{satisfying} \;
\sum_{k=1}^{\infty}b_{k}^{2}\lambda_{k}^{1/2}<+\infty\right\},
\end{split}
\end{equation*}
where $d(x)=\text{\rm dist}(x, \partial \Omega)$, and
$\{\lambda_{k},\varphi_{k} \}$ is the Dirichlet spectral decomposition
of $-\Delta$ in $\Omega$ as above, with $\{\varphi_{k}\}$ an
orthonormal basis of $L^{2}(\Omega)$.

Furthermore, $\V_{0}(\Omega)$ equipped with the norm
\begin{equation}\label{sec1-normV0}
\|u\|_{\V_{0}(\Omega)}=\left\{\|u\|_{H^{1/2}(\Omega)}^{2}
+\int_{\Omega}\frac{u^{2}}{d}\right\}^{1/2}
\end{equation}
is a Banach space.
\end{prop}

The fact that $d^{-1/2}u\in L^2(\Omega)$ if $u$ is the trace of a function in $H_{0,L}^{1}(\C)$ 
follows from a trace boundary Hardy  inequality, originally due to Nekvinda \cite{Nekvi93};
see Lemma~\ref{lem-hardy-inequality} in next section for a simple proof.
Thus, in next section we need to consider 
the operator $A_{1/2}$ defined as in \eqref{defAfirst} but now mapping
$A_{1/2}: \V_{0}(\Omega)\rightarrow \V_{0}^{*}(\Omega)$, where $\V_{0}^{*}(\Omega)$ is the dual space of $\V_{0}(\Omega)$. For $u=\sum_{k=1}^{\infty}b_{k}\varphi_{k}\in \V_{0}(\Omega)$,
we will have $A_{1/2} (\sum_{k=1}^{\infty}b_{k}
\varphi_{k})=\sum_{k=1}^{\infty}b_{k}\lambda_{k}^{1/2}
 \varphi_{k}$
Moreover, there will be a unique harmonic
extension $v\in H_{0,L}^{1}(\C)$ in $\C$ of $u$, and it is given by the expression
\[
v(x,y)=\sum_{k=1}^{\infty}b_{k}\varphi_{k}(x)\exp(-\lambda_{k}^{1/2}y)
\;\;\;\text{for all}\; (x,y)\in \C.
\]
Thus, the operator
$A_{1/2}:\V_{0}(\Omega)\rightarrow \V_{0}^{*}(\Omega)$ is given by
the Dirichlet-Neumann map
\begin{equation*}
 A_{1/2}u:=\frac{\partial v}{\partial \nu}\mid_{\Omega\times\{0\}}
 =\sum_{k=1}^{\infty}b_{k}\lambda_{k}^{1/2}
 \varphi_{k}.
 \end{equation*}
Note that $A_{1/2}\circ A_{1/2}$ is equal to $-\Delta$ in $\Omega$
 with zero Dirichlet boundary value on $\partial \Omega$.
More precisely, we will have that the inverse $B_{1/2}=A_{1/2}^{-1}$ ---which maps $\V_{0}^{*}(\Omega)$ into itself, and also $L^2(\Omega)$ into itself--- is the unique
square root of the inverse Laplacian $(-\Delta)^{-1}$ in $\Omega$
with zero Dirichlet boundary values on $\partial \Omega$; see next section for
details.

 To establish the regularity of weak solutions to (\ref{eqn-frac}) obtained by
 the previous minimization technique,  we establish the following results of
Calder{\'o}n-Zygmund and  of Schauder type  for the linear
problem
\begin{equation}\label{eqn-A-linear0}
 \left\{
\begin{array}{ll} A_{1/2}u=g(x)\quad \quad &\mbox{in} \;\Omega,\\
u=0  &\mbox{on}\;\partial\Omega;
\end{array}
\right.
\end{equation}
see Proposition \ref{prop-linear-reg} for more details.

\begin{teo}\label{teo-Reg} Let $u\in \V_{0}(\Omega)$ be a weak  solution of {\rm(\ref{eqn-A-linear0})}, where $g\in \V_{0}^{*}(\Omega)$ and $\Omega$ is a $C^{2,\alpha}$ bounded domain in $\R^{n}$, for some $0<\alpha<1$.

  If $g\in L^{2}(\Omega)$, then $u\in
H^{1}_{0}(\Omega)$.

 If $g\in H^{1}_{0}(\Omega)$, then $u\in
H^{2}(\Omega)\cap H^{1}_{0}(\Omega)$.

If $g\in L^{\infty}(\Omega)$, then
$u\in C^{\alpha}(\overline{\Omega})$.

If  $g\in
C^{\alpha}(\overline{\Omega})$ and $g|_{\partial \Omega}\equiv 0$,
then $u\in
  C^{1,\alpha}(\overline{\Omega})$.

  If  $g\in C^{1,\alpha}(\overline{\Omega})$ and $g|_{\partial \Omega}\equiv 0$,
  then  $u\in
  C^{2,\alpha}(\overline{\Omega})$.
\end{teo}

In this paper we will give full ---and rather simple--- 
proofs of these regularity results, specially since we could
only find references for some of them and, besides, in close statements to ours 
but not precisely ours. 
Our proof of Theorem
\ref{teo-Reg} uses the extension problem in $\Omega\times (0,\infty)$
related to (\ref{eqn-A-linear0}), and
transforms it to a problem with zero Dirichlet boundary in $\Omega\times\{0\}$
by using an auxiliary function introduced in \cite{CS05}. Then, by
making certain reflections and using classical 
interior regularity theory for the Laplacian, we prove
H\"{o}lder  regularity for $u$ and its derivatives.

To apply the previous H\"{o}lder  regularity linear results to our nonlinear
problem (\ref{eqn-frac}), we first need to prove  that
$g:=f(u)$ is bounded, i.e., $u$ is bounded. We will see that boundedness of weak solutions holds for subcritical and
critical nonlinearities; we establish this result in section~\ref{sec-h-sub}.
We will follow the Brezis-Kato approach bootstrap method. In this way, we establish 
the following
(see Theorem \ref{prop-teo-bk}).

\begin{teo}\label{teo-BK}
Assume that
 $g_{0}$ is a Carath\'{e}odory function in $\Omega\times \R$ satisfying
\[
|g_{0}(x, s)|\le C(1+|s|^{p})\quad\text{for all}\; (x,s)\in
\Omega\times\R,
\]
for some constant $C$,  $1\le p\le \frac{n+1}{n-1}$ if $n\ge 2$, or  $1\le p<\infty$ if
$n=1$,
 where $\Omega$ is a smooth bounded domain in $\R^{n}$.
 Let   $u\in\V_{0}(\Omega)$ be a weak solution of
 \[
\begin{cases}
A_{1/2}u=g_{0}(x,u)&\text{in}\;\Omega, \\u=0
&\text{on}\;\partial\Omega.
\end{cases}
\]
  Then,   $u\in L^{\infty}(\Omega)$.
\end{teo}

The paper is organized as follows. In section \ref{sec-h-pre}, we
study the appropriate function spaces $H_{0,L}^{1}(\C)$ and $\V_{0}(\Omega)$, 
and we give the proof of
Proposition \ref{prop-secone-V0} and other related results. The
regularity results of Theorem \ref{teo-Reg} can be founded in section
\ref{sec-regularity}. Maximum principles, Hopf's lemma, and the
maximum principle in ``small''  domains of Proposition
\ref{prop-smalld0} are proved in section \ref{sec-h-max}. The
complete proof of Theorem \ref{teo-sub1} is given in section
\ref{sec-h-sub} by studying the minimization problem and applying
the previous results on regularity and maximum principles. We prove Theorem
\ref{teo-BK} also in section \ref{sec-h-sub}, while Theorems \ref{teo-ap} and
\ref{thm-l-d2} are established
 in section \ref{sec-h-apriori}, and Theorem
\ref{teo-sys} in section \ref{sec-h-sys}.

\setcounter{equation}{0}
\section[Preliminaries]{Preliminaries: function spaces and the operator $A_{1/2}$}
\label{sec-h-pre}

 In this section we collect preliminary facts for future reference. First of
all, let us set the standard notations to be used in the
paper. We denote the upper half-space in $\R^{n+1}$ by
\[
\R^{n+1}_{+}=\{z=(x,y)=(x_{1},\cdots,x_{n},y)\in \R^{n+1}\mid y>0\}.
\]
Denote by $H^{s}(U)=W^{s,\,2}(U)$ the Sobolev space in a domain $U$
of $\R^{n}$ or of  $\R^{n+1}_{+}$. Letting $U\subset \R^{n}$ and
$s>0$, $H^{s}(U)$ is a Banach space with the norm
\[
\|u\|_{H^{s}(U)}=\Big\{\int_{U}\int_{U}
\frac{|u(x)-u(\bar{x})|^{2}}{|x-\bar{x}|^{n+2s}}\,dxd\bar{x}+\int_{U}|u(x)|^{2}\,dx\Big\}^{1/2}.
\]

Let $\Omega$ be a bounded smooth domain in $\R^{n}$.
 Denote the half-cylinder with base $\Omega$ by
 \[
{\mathcal C}=\Omega\times (0,\infty)
 \]
 and its lateral boundary by
 \[
\partial_{L}{\mathcal C}=\partial\Omega\times [0,\infty).
 \]

To treat the nonlocal problem (\ref{eqn-frac}), we will study a
corresponding extension problem in one more dimension, which allows
us to investigate (\ref{eqn-frac}) by studying a local
problem via classical nonlinear variational methods. We consider the
Sobolev space of functions in $H^{1}(\C)$ whose traces vanish on
$\partial_{L}{\mathcal C}$:
\begin{equation}\label{space}
H_{0,L}^{1}(\C)=\{v\in {H}^{1}({\mathcal{C}})\mid v=0
\;\mbox{a.e.}\;\;\mbox{on}\;
\partial_{L}{\mathcal{C}}\,\},
\end{equation}
equipped with the norm
\begin{equation}\label{norm}
\|v\|=\left(\int_{\C}|\nabla v|^{2}\,dxdy\right)^{1/2}.
\end{equation}
We  denote by $\text{tr}_{\Omega}$ the trace operator on
$\Omega\times\{0\}$ for functions in $H_{0,L}^{1}(\C)$:
\[
\text{tr}_{\Omega}v:=v(\cdot,0),\;\text{for}\; v\in H_{0,L}^{1}(\C).
\]
We have that $\text{tr}_{\Omega}v\in H^{1/2}(\Omega)$, since it is
well known that traces of $H^{1}$ functions are $H^{1/2}$ functions
on the boundary.

Recall the well known spectral theory of the Laplacian $-\Delta$ in
a smooth bounded domain $\Omega$ with zero Dirichlet boundary
values. We repeat each eigenvalue  of $-\Delta$ in $\Omega$ with
zero Dirichlet boundary conditions according to its (finite)
multiplicity:
\[
0<\lambda_{1}<\lambda_{2}\le \cdots\le \lambda_{k}\le
\cdots\rightarrow \infty, \quad \mbox{as}\;k\rightarrow \infty,
\]
 and we denote by $\varphi_{k}\in
H^{1}_{0}(\Omega)$ an eigenfunction corresponding to $\lambda_{k}$
for $k=1,2,\cdots$. Namely,
\begin{equation}\label{eqn-Dirichlet-egv}
\left\{
\begin{array}{ll}
-\Delta \varphi_{k}=\lambda_{k}\varphi_{k}& \mbox{in} \;\Omega,\\
\varphi_{k}=0&\mbox{on}\; \Omega.
\end{array}
\right.
\end{equation}
We can take them to form an orthonormal basis $\{\varphi_{k}\}$ of
$L^{2}(\Omega)$, in particular,
\[
\int_{\Omega}\varphi_{k}^{2}\,dx=1,
\] and to belong to $C^{2}(\overline{\Omega})$ by regularity theory.

Now we can state the main results which we  prove in this section.

\begin{prop}\label{prop-V-def}
Let $\V_{0}(\Omega)$ be the space of all traces on
$\Omega\times\{0\}$ of functions in $H_{0,L}^{1}(\C)$. Then, we
have
\begin{equation*}
\begin{split}
\V_{0}(\Omega)& :=\left\{u =\mbox{{\rm tr}}_{\Omega}v\mid v\in
H_{0,L}^{1}(\C)\right\}
\\
&=\left\{u\in H^{1/2}(\Omega)\mid\int_{\Omega}\frac{u^{2}(x)}{d(x)}\,dx<+\infty\right\}\\
&= \left\{ u \in L^{2}(\Omega) \mid
u=\sum_{k=1}^{\infty}b_{k}\varphi_{k} \;\text{satisfying} \;
\sum_{k=1}^{\infty}b_{k}^{2}\lambda_{k}^{1/2}<+\infty\right\},
\end{split}
\end{equation*}
where $d(x)=\text{\rm dist}(x, \partial \Omega)$, and
$\{\lambda_{k},\varphi_{k} \}$ is the Dirichlet spectral decomposition
of $-\Delta$ in $\Omega$ as above, with $\{\varphi_{k}\}$ an
orthonormal basis of $L^{2}(\Omega)$.

Furthermore, $\V_{0}(\Omega)$ equipped with the norm
\begin{equation}\label{normV0}
\|u\|_{\V_{0}(\Omega)}=\left\{\|u\|_{H^{1/2}(\Omega)}^{2}
+\int_{\Omega}\frac{u^{2}}{d}\right\}^{1/2}
\end{equation}
is a Banach space.
\end{prop}

\begin{prop}\label{prop-A-half-def}
If $u\in \V_{0}(\Omega)$, then there exists a unique harmonic
extension $v$ in $\C$ of $u$ such that $v\in H_{0,L}^{1}(\C)$. In
particular, if the expansion of $u$ is written by
$u=\sum_{k=1}^{\infty}b_{k}\varphi_{k}\in \V_{0}(\Omega)$, then
\[
v(x,y)=\sum_{k=1}^{\infty}b_{k}\varphi_{k}(x)\exp(-\lambda_{k}^{1/2}y)\;\;\;\text{for
all}\; (x,y)\in \C,
\]
 where
$\{\lambda_{k},\varphi_{k} \}$ is the Dirichlet spectral decomposition
of $-\Delta$ in $\Omega$ as above, with $\{\varphi_{k}\}$ an
orthonormal basis of $L^{2}(\Omega)$.

The operator $A_{1/2}:\V_{0}(\Omega)\rightarrow \V_{0}^{*}(\Omega)$
is given by
\begin{equation*}
 A_{1/2}u:=\left.\frac{\partial v}{\partial \nu}\right|_{\Omega\times\{0\}},
 \end{equation*}
 where $\V_{0}^{*}(\Omega)$ is the dual space of $\V_{0}(\Omega)$.
We have that
\[
A_{1/2}u=\sum_{k=1}^{\infty}b_{k}\lambda_{k}^{1/2}
 \varphi_{k},\]
 and that $A_{1/2}\circ A_{1/2}$ {\rm(}when $A_{1/2}$ is acting, for instance,
 on smooth functions with compact support in $\Omega${\rm)}
  is equal to $-\Delta$ in $\Omega$
 with zero Dirichlet boundary values on $\partial \Omega$.
More precisely, the inverse $B_{1/2}:=A_{1/2}^{-1}$ is the unique positive
square root of the inverse Laplacian $(-\Delta)^{-1}$ in $\Omega$
with zero Dirichlet boundary values on $\partial \Omega$.
\end{prop}

The proofs of these two propositions need the development of several
tools.   First let us give some properties of the space
$H_{0,L}^{1}(\C)$. Denote by ${\mathcal{D}}^{1,2}(\R^{n+1}_{+})$
the closure of the set of  smooth functions compactly supported in
$\overline{\R^{n+1}_{+}}$ with respect to the norm of
$\|w\|_{{\mathcal{D}}^{1,2}(\R^{n+1}_{+})}=\Big(\int_{\R^{n+1}_{+}}|\nabla
w|^{2}\,dxdy\Big)^{1/2}$.  We recall the well known Sobolev trace inequality
that for $w\in {\mathcal{D}}^{1,\,
2}(\R^{n+1}_{+})$,
\begin{equation}\label{eqn-inequ-ST}
\Big(\int_{\R^{n}}|w(x,0)|^{2n/(n-1)}dx\Big)^{(n-1)/2n}\le
C\Big(\int_{\R_{+}^{n+1}}|\nabla w(x,y)|^{2}dxdy\Big)^{1/2},
\end{equation}
where $C$ depends only on $n$.

Denote for $n\ge 2$,
\[ 2^{\sharp}=\frac{2n}{n-1}\quad\text{and}\quad 2^{\sharp}-1=\frac{n+1}{n-1}.
\]
We say that $p$ is subcritical if $1< p<
2^{\sharp}-1=\frac{n+1}{n-1}$ for $n\ge 2$, and $1< p<\infty$ for
$n=1$. We also say that $p$ is critical if $p=
2^{\sharp}-1=\frac{n+1}{n-1}$ for $n\ge 2$, and that $p$ is
supercritical if $p> 2^{\sharp}-1=\frac{n+1}{n-1}$ for $n\ge 2$.

Lions \cite{Lions85} showed that
\begin{equation}\label{eqn-trace-S}
S_{0}=\inf\left\{ \frac{\int_{\R^{n+1}_{+}}|\nabla
w(x,y)|^{2}dxdy}{(\int_{\R^{n}}|w(x,0)|^{2^{\sharp}}dx)^{2/2^{\sharp}}}
\mid w\in {\mathcal{D}}^{1,2}(\R^{n+1}_{+})\, \right\}
\end{equation}
is achieved.
Escobar \cite{Es} prove that the extremal functions have all the
form
\begin{equation}\label{eqn-U}
U_{\varepsilon}(x,y)=\frac{\varepsilon^{(n-1)/2}}{|(x-x_{0},y+\varepsilon)|^{n-1}},
\end{equation}
where $x_{0}\in \R^{n}$ and $\varepsilon >0$ are arbitrary. In addition, the best constant is
\[
S_{0}=\frac{(n-1)\sigma_{n}^{1/n}}{2},
\]
where $\sigma_{n}$ denotes the volume of $n$-dimensional sphere ${\mathbb{S}}^{n}\subset \R^{n+1}$.

The Sobolev trace inequality leads directly to the next three lemmas.
For $v\in H_{0,L}^{1}(\C)$, its extension by zero in
$\R^{n+1}_{+}\setminus\C$ can be approximated by functions compactly
supported in $\overline{\R^{n+1}_{+}}$. Thus the
 Sobolev trace inequality (\ref{eqn-inequ-ST}) leads to:
\begin{lema}\label{pstar-in} Let $n\ge 2$ and
$2^{\sharp}=\frac{2n}{n-1}$. Then there exists a constant $C$,
depending only on $n$, such that, for all $v\in H_{0,L}^{1}(\C)$,
\begin{equation}\label{eqn-in-sob-trace}
\Big(\int_{\Omega}|v(x,0)|^{2^{\sharp}}dx\Big)^{1/2^{\sharp}}\le
C\Big(\int_{{\mathcal{C}}}|\nabla v(x,y)|^{2}dxdy\Big)^{1/2}.
\end{equation}
\end{lema}

 By H\"{o}lder's inequality, since $\Omega$ is bounded,
the above lemma leads to:
\begin{lema}\label{p-in} Let $1\le q \le 2^{\sharp}$ for $n\ge 2$.
Then, we have that for all $v\in H_{0,L}^{1}(\C)$,
\begin{equation}\label{eqn-h-inq-2}
\Big(\int_{\Omega}|v(x,0)|^{q}dx\Big)^{1/q}\le
C\Big(\int_{{\mathcal{C}}}|\nabla v(x,y)|^{2}dxdy\Big)^{1/2},
\end{equation}
 where $C$ depends only on $n, q$, and the measure of $\Omega$.
Moreover, {\rm(\ref{eqn-h-inq-2})} also holds for
  $1\le q <\infty$ if  $n=1$.
\end{lema}

This lemma states that $\text{tr}_{\Omega}(H_{0,L}^{1}(\C))\subset
L^{q}(\Omega)$, where $1\le q \le 2^{\sharp}$ for $n\ge 2$ and $1\le
q<\infty$ for $n=1$,  (see the proof of Lemma \ref{prop-embedding} for the
case $n=1$). In addition, we also have the following compact
embedding.
\begin{lema}\label{prop-embedding} Let $1\le q < 2^{\sharp}=\frac{2n}{n-1}$
for $n\ge 2$ and $1\le q < \infty$ for $n=1$.
 Then
 $\text{\rm{tr}}_{\Omega}(H_{0,\,L}^{1}(\C))$ is compactly embedded in
$L^{q}(\Omega)$.
\end{lema}
\noindent{\bf Proof.}  It is well known that
$\text{tr}_{\Omega}(H_{0,L}^{1}(\C))\subset H^{1/2}(\Omega)$ and
that $H^{1/2}(\Omega)$ $\subset\subset L^{q}(\Omega)$ when $1\le q <
2^{\sharp}=\frac{2n}{n-1}$ for $n\ge 2$ and $1\le q < \infty$ for
$n=1$. Here $\subset\subset$ denotes the compact embedding. This
completes the proof of the lemma. However, if one wants to avoid the
use of the fractional Sobolev space $H^{1/2}(\Omega)$, the following
is an alternative simple proof.

Considering the restriction of functions in $\C$ to
$\Omega\times(0,1)$, it suffices to show that the embedding is
compact with $\C$ replaced by $\Omega\times(0,1)$. To prove this,
let $v_{m}\in H^{1}_{0,L}(\Omega\times (0,1)):=\{v\in
H^{1}(\Omega\times (0,1))\mid v=~0\; \text{a.e.}\; \text{on}\;
\partial{\Omega}\times(0,1)\}$ such that
$v_{m}\rightharpoonup 0 $ weakly in $ H^{1}_{0,L}(\Omega\times
(0,1))$, as $m\rightarrow \infty$.
 We may
assume by the classical Rellich's theorem in $\Omega\times(0,1)$
that $v_{m} \rightarrow 0 $ strongly in $ L^{2}(\Omega\times
(0,1))$, as $m\rightarrow \infty$. We introduce the function
$w_{m}=(1-y)v_{m}$. It is clear that
\[ w_{m}|_{\Omega\times\{0\}}=v_{m},\quad w_{m}|_{\Omega\times\{1\}}=0.
\]
By direct computations we have
\begin{equation*}
\begin{split}
&\int_{\Omega}|v_{m}(x,0)|^{2}\,dx=\int_{\Omega}|w_{m}(x,0)|^{2}\,dx
=-\int_{0}^{1}\int_{\Omega}\partial_{y}(w_{m}^{2}(x,y))\,dxdy\\
&\le
2\Big(\int_{0}^{1}\int_{\Omega}w_{m}^{2}(x,y)\,dxdy\Big)^{1/2}\Big(\int_{0}^{1}\int_{\Omega}|\nabla
w_{m}(x,y)|^{2}\,dxdy\Big)^{1/2}.
 \end{split}
\end{equation*}
Therefore, since  $w_{m}=(1-y)v_{m}$ is bounded in
$H^{1}(\Omega\times (0,1))$ and $w_{m}\rightarrow 0$ strongly in $L^{2}(\Omega\times(0,1))$, we find that, as $m\rightarrow \infty$,
\begin{equation*}
v_{m}(x,0)\rightarrow 0\quad \mbox{ strongly in}\;
L^{2}(\Omega)\;\mbox{and hence also in} \; L^{1}(\Omega).
\end{equation*}
On the other hand, since $q$ is subcritical, the following
interpolation inequality,
\[
\|v_{m}(\cdot,0)\|_{L^{q}(\Omega)}\le
\|v_{m}(\cdot,0)\|^{\theta}_{L^{1}(\Omega)}\|v_{m}(\cdot,0)\|_{L^{2^{\sharp}}(\Omega)}^{1-\theta}
\]
for some $0<\theta<1$ completes the proof since we already know that
$v_{m}$  converges strongly to zero in $L^{1}(\Omega)$.
\hfill{$\Box$}

\medskip

We also need to establish a trace boundary Hardy inequality, which already appeared in a work of
Nekvinda \cite{Nekvi93}.

\begin{lema}\label{lem-hardy-inequality} We have that
\[
\text{{\rm tr}}_{\Omega} (H_{0,L}^{1}(\C))\subset H^{1/2}(\Omega)
\] is a continuous injection. In addition, for every  $v\in H_{0,L}^{1}(\C)$,
\begin{equation}\label{eqn-bh-ineq}
\int_{\Omega} \frac{|v(x,0)|^{2}}{d (x)}\,dx \le C\int_{\C}|\nabla
v(x,y)|^{2}\,dxdy,
\end{equation}
where $d(x)=\mbox{{\rm dist}}
(x,\partial \Omega)$ and the constant $C$ depends only on $\Omega$.
\end{lema}
\noindent{\bf Proof.} The first statement is clear since the traces
of $H^{1}(\C)$ functions belong to $H^{1/2}(\partial \C)$. Regarding
the second statement, we prove it in two steps.\\
 {\it Step 1.} Assume first that $n=1$ and
$\Omega=(0,1)$.   For $0<x_{0}<1/2$, consider the
 segment from $(0,x_{0})$ to $(x_{0},0)$ in $\C=(0,1)\times
(0,\infty)$. We have
\[
v(x_{0},0)=v(t,x_{0}-t)\mid_{t=0}^{x_{0}}=\int_{0}^{x_{0}}
(\partial_{x}v-\partial_{y}v)(t,x_{0}-t)\,dt.
\]
Then
\[
|v(x_{0},0)|^{2}\le x_{0}\int_{0}^{x_{0}}2|\nabla
v(t,x_{0}-t)|^{2}\,dt.
\]
Dividing this inequality by $x_{0}$ and integrating in $x_{0}$ over
$(0,1/2)$, and making the change of variables $x=t$, $y=x_{0}-t$, we
deduce
\[
\int_{0}^{1/2}\frac{|v(x_{0},0)|^{2}}{x_{0}}\, dx_{0} \le
2\int_{0}^{1/2}dx\int_{0}^{1/2}dy|\nabla v|^{2} \le
2\int_{\C}|\nabla v|^{2}\,dxdy.
\]
Doing the same on $(1/2,1)$, this establishes inequality
(\ref{eqn-bh-ineq}) of the lemma.\\
 {\it Step 2.} In the general case, after straightening  a piece of the
boundary $\partial\Omega$ and rescaling the new variables, we can
consider the inequality in a domain $D=\{x=(x',x_{n})\mid |x'|<1,
0<x_{n}<1/2\}$ and assume that $v=0$ on
$\{x_{n}=0,|x'|<1\}\times(0,\infty)$, since the flatting procedure
possesses equivalent norms. By the argument in  Step 1 above,
we have
\[
\int_{0}^{1/2}\frac{|v(x,0)|^{2}}{x_{n}}\, dx_{n} \le
C\int_{0}^{1/2}\int_{0}^{\infty}|\nabla v|^{2}\,dx_{n}dy,
\]
for all $x'$ with $|x'|<1$. From this, integrating in $x'$ we have
\begin{align*}
\int_{D}\frac{|v(x,0)|^{2}}{x_{n}}\,dx
=&\int_{D}\int_{0}^{1/2}\frac{|v(x,0)|^{2}}{x_{n}}\,dx'dx_{n}\\
\le& C\int_{D\times(0,\infty)}|\nabla v|^{2}\,dxdy.
\end{align*}
Since after flattening of $\partial \Omega$, $x_{n}$ is comparable
to $d(x)=\text{dist}(x,\partial \Omega)$, this is the desired
inequality (\ref{eqn-bh-ineq}). \hfill \hfill{$\Box$}

\medskip

Recall that the fractional Sobolev space $H^{1/2}(\Omega)$ is a
Banach space with the norm
\begin{equation}\label{eqn-sobolevhalf-norm}
\|u\|^{2}_{H^{1/2}(\Omega)}=\int_{\Omega}\int_{\Omega}
\frac{|u(x)-u(\bar{x})|^{2}}{|x-\bar{x}|^{n+1}}\,dxd\bar{x}+\int_{\Omega}|u(x)|^{2}\,dx.
\end{equation}
 Note that the closure $H_{0}^{1/2}(\Omega)$ of smooth functions with compact support,
  $C_{c}^{\infty}(\Omega)$, in $H^{1/2}(\Omega)$
 is all the space $H^{1/2}(\Omega)$, (see Theorem 11.1 in \cite{LM}). That is,
 $C_{c}^{\infty}(\Omega)$ is dense in  $H^{1/2}(\Omega)$. However, in contrast with this, the trace in
 $\Omega$ of functions in $H_{0,L}^{1}(\C)$ ``vanish'' on $\partial \Omega$ in the sense  given by
 (\ref{eqn-bh-ineq}).

 Recall that we have denoted  by $\V_{0}(\Omega)$
 the space of traces on $\Omega\times\{0\}$ of functions
 in $H_{0,L}^{1}(\C)$:
\begin{equation}\label{def-V-zero}
\V_{0}(\Omega) :=\{u =\mbox{{\rm tr}}_{\Omega}v\mid v\in
H_{0,L}^{1}(\C)\} \subset H^{1/2}(\Omega),
\end{equation}
endowed with the norm (\ref{normV0}) in Proposition \ref{prop-V-def}. The dual space of
$\V_{0}(\Omega)$ is denoted by $\V_{0}^{*}(\Omega)$, equipped  with
the norm
\[
\|g\|_{\V_{0}^{*}(\Omega)}=\sup\{\langle u,g \rangle\mid u\in
\V_{0}(\Omega), \|u\|_{\V_{0}(\Omega)}\le 1\}.
\]

Next we give the first characterization  of the space
$\V_{0}(\Omega)$:
\begin{lema}\label{lem-V-zero-withd} Let $\V_{0}(\Omega)$
be the space of traces on $\Omega\times\{0\}$ of functions in
$H_{0,L}^{1}(\C)$, as in {\rm(\ref{def-V-zero})}. Then, we
have
\[
\V_{0}(\Omega)=\left\{u\in
H^{1/2}(\Omega)\mid\int_{\Omega}\frac{u^{2}(x)}{d(x)}\,dx<+\infty\right\},
\]
where $d(x)=\mbox{{\rm dist}} (x,\partial \Omega)$.
\end{lema}
\noindent{\bf Proof.} The inclusion $\subset$ follows Lemma
\ref{lem-hardy-inequality}. Next we show the other inclusion.
Let $u\in H^{1/2}(\Omega)$ satisfy $\int_{\Omega}u^{2}/d<\infty$.
 Let $\tilde{u}$ be the extension of $u$ in all of
$\R^{n}$ assigning $\tilde{u}\equiv 0$ in $\R^{n}\setminus \Omega$. The quantity
\[
\|\tilde{u}\|^{2}_{H^{1/2}(\R^{n})}=\int_{\R^{n}}\int_{\R^{n}}
\frac{|\tilde{u}(x)-\tilde{u}(\bar{x})|^{2}}{|x-\bar{x}|^{n+1}}\,dxd\bar{x}+\int_{\R^{n}}|\tilde{u}(x)|^{2}\,dx
\]
 can be bounded ---using
$\tilde{u}\equiv 0$ in $\R^{n}\setminus \Omega$--- by a constant
times
\[\left\{\|u\|^{2}_{H^{1/2}(\Omega)}+\int_{\Omega}\frac{u^{2}(x)}{d(x)}\,dx\right\}^{1/2},\]
 that we assume to be finite. Hence, $\tilde{u}\in H^{1/2}(\R^{n})$
and thus $\tilde{u}$ is the trace in $\R^{n}=\partial \R^{n+1}_{+}$
of a function $\tilde{v}\in H^{1}(\R^{n+1}_{+})$.

Next, we use a partition of the unity, and local bi-Lipschitz maps (defined below) sending
$\overline{\R^{n+1}_{+}}$ into $\overline{\Omega}\times
[0,\infty)=\overline{\C}$ being the identity on
$\Omega\times\{0\}$  and mapping  $\R^{n}\setminus \Omega =
(\partial \R^{n+1}_{+})\setminus \Omega$ into $\partial \Omega
\times [0,\infty)$. By composing these maps with the function (cutted off with the partition of unity) $\tilde{v}$, we
obtain an $H_{0,L}^{1}(\C)$ function with $u$ as trace on
$\Omega\times\{0\}$, as desired.

Finally we give a concrete expression for one such bi-Lipschitz
maps. First,
 consider the one dimensional case $\Omega=(0,\infty)$. Then simply take the bi-Lipschitz map
\begin{equation*}
\begin{split}
 (x,y)&\in (0,\infty)\times (0,\infty)=\Omega\times (0,\infty)\\
&\mapsto (\frac{x^{2}-y^{2}}{\sqrt{x^{2}+y^{2}}},
\frac{2xy}{\sqrt{x^{2}+y^{2}}}) \in \R\times(0,\infty),
\end{split}
\end{equation*}
whose Jacobian can be checked to be identically $2$.
 In the general case, we can flatten the boundary $\partial\Omega$
  and use locally the previous map.
 \hfill{$\Box$}

\medskip

Next we consider, for a given function $u\in \V_{0}(\Omega)$, the minimizing problem:
\begin{equation}\label{lem-min-ext}
\inf\left\{\int_{\C}|\nabla v|^{2}\,dxdy\mid v\in H_{0,L}^{1}(\C),
v(\cdot,0)=u \;\mbox{in}\; \Omega\right\}.
\end{equation}
By the definition of $\V_{0}(\Omega)$, the set of functions $v$ where we minimize is non empty. By lower weak semi-continuity and by Lemma
\ref{prop-embedding}, we see that there exists a minimizer $v$. We
will prove next that this minimizer $v$ is unique.
 We call $v$  a {\it weak solution} of the problem
\begin{equation}\label{eqn-extension}
\left\{
\begin{array}{ll}
 \Delta v=0 &  \mbox{in}\; \C, \\
 v=0 & \mbox{on}\; \partial_{L}{\mathcal C},\\
v=u & \mbox{on} \; \Omega\times \{0\}.
\end{array}
\right.
\end{equation}
That is, we have
\begin{lema}\label{lem-extension1}
For $u\in \V_{0}(\Omega)$, there exists a unique minimizer $v$ of
{\rm(\ref{lem-min-ext})}. The function $v\in H^{1}_{0,L}(\C)$ is
the harmonic extension of $u$ (in the weak sense) to $\C$ vanishing
on $\partial_{L} \C$.
\end{lema}
\noindent{\bf Proof.} By the definition of $\V_{0}(\Omega)$, we have
that, for every $u\in \V_{0}(\Omega)$,  there exists at least one
$w\in H_{0,L}^{1}(\C)$ such that $\mbox{{\rm tr}}_{\Omega}(w)=u$.
Then the standard minimization argument gives  (using lower
semi-continuity and Lemma \ref{prop-embedding}) the existence of a
minimizer. The uniqueness of minimizer follows automatically from the identity of
the parallelogram used
 for two possible minimizers $v_{1}$ and
$v_{2}$,
\[
0 \le
J(\frac{v_{1}-v_{2}}{2})=\frac{1}{2}J(v_{1})+\frac{1}{2}J(v_{2})-J(\frac{v_{1}+v_{2}}{2})\le
0,
\]
where $J(v)=\int_{\C}|\nabla v|^{2}\,dxdy$, which leads to $v_{1}=v_{2}$.
 \hfill{$\Box$}

\medskip

By Lemma \ref{lem-extension1}, there exists a unique function $v\in
H_{0,L}^{1}(\C)$ which is the harmonic extension of $u$ in $\C$
vanishing on $\partial_{L}\C$, and that we denote by
\[
v:=\mbox{h-ext}(u).
\]
It is easy to see that for every $\eta\in C^{\infty}(\overline{\C})\cap H^{1}(\C)$
and $\eta\equiv 0$ on $\partial_{L}\C$,
\begin{equation}\label{eqn-Xi-V-weak}
\int_{\C}\nabla v\nabla \eta\,dxdy=\int_{\Omega}\frac{\partial
v}{\partial \nu}\eta\,dx.
\end{equation}
 By Lemma \ref{lem-hardy-inequality}, there exists a constant $C$ such that
 for every $u\in \V_{0}(\Omega)$,
\begin{equation}\label{eqn-bi-one1}
 \|u\|_{\V_{0}(\Omega)} \le C\|\mbox{h-ext}(u)\|_{H_{0,L}^{1}(\C)}.
\end{equation}
Next, note that the h-ext operator is bijective from $\V_{0}(\Omega)$ to the subspace $\mathcal{H}$ of
$H_{0,L}^{1}(\C)$ formed by all harmonic functions in $H_{0,L}^{1}(\C)$. Since
both $\V_{0}(\Omega)$ and $\mathcal{H}$ are Banach spaces, the open mapping theorem  gives that we also have the reverse inequality to (\ref{eqn-bi-one1}), i.e., there exists a constant $C$ such that
\begin{equation}\label{eqn-bi-one2}
\|\text{h-ext}(u)\|_{H_{0,L}^{1}(\C)}\le
C\|u\|_{\V_{0}(\Omega)},
\end{equation}
for all $u\in \V_{0}(\Omega)$. From this we deduce the following.
Given a smooth $\xi\in \V_{0}(\Omega)$, consider the
$\text{h-ext}(\xi)$ and call it $\eta$. Now, we use
(\ref{eqn-Xi-V-weak}) and (\ref{eqn-bi-one2}) (for $u$ and $\xi$) to obtain $\Big|\int_{\Omega}\frac{\partial
v}{\partial \nu}\xi\,dx\Big|\le
C\|u\|_{\V_{0}(\Omega)}\|\xi\|_{\V_{0}(\Omega)}$. That is,
   $\frac{\partial v}{\partial \nu}\mid_{\Omega}\in \V_{0}^{*}(\Omega)$
   and there is the bound:
   \[
   \left\|\frac{\partial}{\partial \nu} \text{h-ext}(u)\right\|_{\V_{0}^{*}(\Omega)}\le C\|u\|_{\V_{0}(\Omega)}.
   \]
   Hence we have
\begin{lema}
The operator $A_{1/2}: \V_{0}(\Omega)\rightarrow\V_{0}^{*}(\Omega)$
defined by
 \begin{equation}\label{def-A}
 A_{1/2}u:=\left.\frac{\partial v}{\partial \nu}\right|_{\Omega\times\{0\}},
 \end{equation}
where $v=\mbox{\rm h-ext}(u)\in H_{0,L}^{1}(\C)$ is the harmonic
extension of $u$ in $\C$ vanishing on $\partial_{L}\C$,  is linear and bounded from $\V_{0}(\Omega)$ to
$\V_{0}^{*}(\Omega)$.
\end{lema}

We now give the spectral representation of $A_{1/2}$ and the
corresponding structure of the space $\V_{0}(\Omega)$.
\begin{lema}\label{lem-spectrum} (i) Let $\{\varphi_{k}\}$
be an orthonormal basis  of $L^{2}(\Omega)$ forming a spectral
decomposition of $-\Delta$ in $\Omega$ with Dirichlet boundary
conditions as in {\rm(\ref{eqn-Dirichlet-egv})},
with $\{\lambda_{k}\}$ the corresponding Dirichlet eigenvalues of
$-\Delta$ in $\Omega$. Then, we have
\[
\V_{0}(\Omega)=\left\{ u=\sum_{k=1}^{\infty}b_{k}\varphi_{k}\in
L^{2}(\Omega)\mid
\sum_{k=1}^{\infty}b_{k}^{2}\lambda_{k}^{1/2}<+\infty \right\}.
\]
(ii) Let $u\in \V_{0}(\Omega)$. Then we have, if
$u=\sum_{k=1}^{\infty}b_{k}\varphi_{k}$,
\[
A_{1/2}u=\sum_{k=1}^{\infty}b_{k}\lambda_{k}^{1/2}
 \varphi_{k}\in \V_{0}^{*}(\Omega).
\]
\end{lema}
\noindent{\bf Proof.}  Let  $u\in \V_{0}(\Omega)$, which is
contained in $L^{2}(\Omega)$. Let its expansion be written by
$u(x)=\sum_{k=1}^{\infty}b_{k}\varphi_{k}(x)$. Consider  the
function
\begin{equation}\label{eqn-f-v}
v(x,y)=\sum_{k=1}^{\infty}b_{k}\varphi_{k}(x)\exp(-\lambda_{k}^{1/2}y),
\end{equation}
which is clearly smooth for $y>0$. Observe that   $v(x,0)=u(x)$ in
$\Omega$ and, for $y>0$,
\[
\Delta
v(x,y)=\sum_{k=1}^{\infty}b_{k}\{-\lambda_{k}\varphi_{k}(x)\exp(-\lambda_{k}^{1/2}y)
+\lambda_{k}\varphi_{k}(x)\exp(-\lambda_{k}^{1/2}y)\}=0.
\]
Thus, $v$ is a harmonic extension of $u$. We will have that $v=\mbox{h-ext}(u)$,
by uniqueness, once we find the condition on $\{b_{k}\}$ for $v$ to
belong to $H_{0,L}^{1}(\C)$. But such condition is simple. Using
(\ref{eqn-f-v}) and that $\{\varphi_{k}\}$ are eigenfunctions of
$-\Delta$ and orthonormal in $L^{2}(\Omega)$, we have
\begin{equation*}
\begin{split}
\int_{0}^{\infty}\int_{\Omega}|\nabla v|^{2}\,dxdy&=
\int_{0}^{\infty}\int_{\Omega}\{|\nabla_{x}v|^{2}+|\partial_{y}v|^{2}\}\,dxdy\\
&=2\sum_{k=1}^{\infty}b_{k}^{2}\lambda_{k}\int_{0}^{\infty}
\exp(-2\lambda_{k}^{1/2}y)\,dy
\\&=2\sum_{k=1}^{\infty}b_{k}^{2}\lambda_{k}\frac{1}{2\lambda_{k}^{1/2}}
=\sum_{k=1}^{\infty}b_{k}^{2}\lambda_{k}^{1/2}.
\end{split}
\end{equation*}
This means that $v\in H^{1}_{0,L}(\C)$ if and only if
$\sum_{k=1}^{\infty}b_{k}^{2}\lambda_{k}^{1/2}<\infty$. Therefore,
this condition on $\{b_{k}\}$ is equivalent to $u\in
\V_{0}(\Omega)$.

Assertion (ii) follows from the direct computation of
$-\frac{\partial v}{\partial y}\mid_{y=0}$ using (\ref{eqn-f-v}).
\hfill{$\Box$}

\medskip 


In functional analysis, the classical spectral decomposition holds for self-adjoint compact operators, such as the Dirichlet inverse Laplacian $(-\Delta)^{-1}:L^{2}(\Omega)\rightarrow L^{2}(\Omega)$. This is the reason why we now define, with the aid  of the Lax-Milgram 
theorem, a compact operator $B_{1/2}$ which will be the inverse of $A_{1/2}$.
\begin{definition}
{\rm
Define the operator $B_{1/2}:
\V^{*}_{0}(\Omega)\rightarrow\V_{0}(\Omega) $, by $g \mapsto
\mbox{{\rm tr}}_{\Omega} v$, where $v$ is found by solving the
problem:
\begin{equation}\label{eqn-B-half}
\left\{
\begin{array}{ll}
\Delta v=0 &  \mbox{in}\; \C, \\
 v=0 & \mbox{on}\; \partial_{L}{\mathcal C},\\
\frac{\partial v}{\partial \nu}=g(x) & \mbox{on} \; \Omega\times
\{0\},
\end{array}
\right.
\end{equation}
as we indicate next.
}
\end{definition}
We say that $v$ is a weak solution of (\ref{eqn-B-half}) whenever
$v\in H_{0,L}^{1}(\C)$ and
\begin{equation}\label{eqn-sol-w}
\int_{\C}\nabla v\nabla \xi\,dxdy=\langle g,\xi(\cdot,0)\rangle
\end{equation}
for all $\xi\in H_{0,L}^{1}(\C)$. We see that there exists  a
 unique weak solution of (\ref{eqn-B-half}) by the Lax-Milgram theorem, via studying the corresponding
 functional in $H_{0,\,L}^{1}(\C)$:
 \[
I(v)=\frac{1}{2}\int_{\C}|\nabla v|^{2}\,dxdy-\langle
g,v(\cdot,0)\rangle,
 \]
 where $g\in \V_{0}^{*}(\Omega)$ is given.
 Observe  that the operator $B_{1/2}$ is clearly the inverse of the operator
 $A_{1/2}$.

On the other hand, let us compute $B_{1/2}\circ B_{1/2}\!\!\mid_{L^{2}(\Omega)}$. Here note that
since $\V_{0}(\Omega)\subset L^{2}(\Omega)$, we have $L^{2}(\Omega)\subset \V_{0}^{*}(\Omega)$. For
a given $g\in L^{2}(\Omega)$, let
$\varphi\in H^{1}_{0}(\Omega)\cap H^{2}(\Omega)$ be the solution of
Poisson's problem for the Laplacian
\[
\left\{
\begin{array}{ll}
 -\Delta \varphi=g &  \mbox{in}\; \Omega, \\
 \varphi=0 & \mbox{on}\; \partial\Omega.
\end{array}
\right.
\]
 Since $H^{1}_{0}(\Omega)\subset
\V_{0}(\Omega)$ (for instance, by Lemma \ref{lem-spectrum}), there is 
a unique  harmonic extension $\psi\in H_{0,L}^{1}(\C)$ of $\varphi$
in $\C$ such that
\[
\left\{
\begin{array}{ll}
 \Delta \psi=0 &  \mbox{in}\; \C, \\
 \psi=0 & \mbox{on}\; \partial_{L}{\mathcal C},\\
\psi=\varphi & \mbox{on} \; \Omega\times \{0\}.
\end{array}
\right.
\]
Moreover, $\tilde{\psi}(x,y):=\psi(x,y)-\varphi(x)$ solves
\[
\left\{
\begin{array}{ll}
 -\Delta \tilde{\psi}=\Delta \varphi=-g(x) &  \mbox{in}\; \C, \\
 \tilde{\psi}=0 & \mbox{on}\;\partial_{L}{\mathcal C},\\
\tilde{\psi}=0& \mbox{on} \;\Omega\times \{0\}.
\end{array}
\right.
\]
Considering the odd reflection $\widetilde{\psi}_{od}$
of $\widetilde{\psi}$  across $\Omega\times\{0\}$,
and the function
\[
g_{od}(x,y)=
\left\{
\begin{array}{ll}
g(x),&y\ge 0,\\
-g(x),&y<0,
\end{array}
\right.
\]
  we have
\[
\left\{
\begin{array}{ll}
 -\Delta \widetilde{\psi}_{od}=-g_{od} &  \mbox{in}\; \Omega\times \R, \\
 \widetilde{\psi}_{od}=0 & \mbox{on}\;\partial\Omega\times \R.
\end{array}
\right.
\]
Therefore, since $g_{od}\in L^{2}(\Omega\times(-2,2))$, we deduce
$\widetilde{\psi}_{od}\in H^{2}(\Omega\times(-1,1))$ and hence
$\psi\in H^{2}(\Omega\times(0,1))$. We deduce, by the smoothness of the harmonic function
$\psi$ for $y>0$ and by its exponential
decay in $y$ ---see (\ref{eqn-f-v})---, that $\psi\in
H_{0,\,L}^{1}(\C)\cap H^{2}(\C)$.

It follows that $-\partial_{y}\psi\in H_{0,L}^{1}(\C)$ solves
 \[
\left\{
\begin{array}{ll}
 \Delta( -\partial_{y}\psi)=0 & \mbox{in}\; \C, \\
-\partial_{y}\psi=0& \mbox{on}\; \partial_{L}{\mathcal C},
\end{array}
\right.
\]
and
\[
\frac{\partial}{\partial
\nu}(-\partial_{y}\psi)=\partial_{yy}\psi=-\Delta_{x}\psi=-\Delta\varphi=g
  \quad \mbox{on} \;\Omega\times \{0\}.
\]
 Since $\V_{0}(\Omega)\subset L^{2}(\Omega)$, we have that
$g\in L^{2}(\Omega)\cong L^{2}(\Omega)^{*}\subset
\V_{0}^{*}(\Omega)$, and we deduce that the solution $v\in \V_{0}(\Omega)$ of
(\ref{eqn-B-half}) is $v=-\partial_{y}\psi$, because of the
uniqueness of $H^{1}_{0,L}(\C)$ solution of
(\ref{eqn-B-half}). In particular,
$B_{1/2}g=v(\cdot,0)=-\partial_{y}\psi(\cdot,0)$. On the other hand,
since $\psi\in H_{0,L}^{1}(\C)$ solves
\[
\left\{
\begin{array}{ll}
 \Delta \psi=0 &  \mbox{in}\; \C, \\
 \psi=0 & \mbox{on}\; \partial_{L}{\mathcal C},\\
\frac{\partial \psi}{\partial \nu}\equiv
-\partial_{y}\psi(\cdot,0)=v(\cdot,0)=B_{1/2}g & \mbox{on} \; \Omega\times
\{0\},
\end{array}
\right.
\]
we conclude that
\[
(B_{1/2}\circ B_{1/2})g=
B_{1/2}v(\cdot,0)=\psi(\cdot,0)=\varphi=(-\Delta)^{-1}g.
\]

Summarizing the above argument, we have:
\begin{prop}\label{prop-B-composition}
 $B_{1/2}\circ B_{1/2}\!\!\mid_{L^{2}(\Omega)}
 =(-\Delta)^{-1}: L^{2}(\Omega)\rightarrow L^{2}(\Omega)$,
 where $(-\Delta)^{-1}$ is
 the inverse Laplacian in $\Omega$ with zero Dirichlet boundary conditions.
\end{prop}

Note that $B_{1/2}:L^{2}(\Omega)\rightarrow L^{2}(\Omega)$ is a
self-adjoint operator. In fact, since for $v_{1},v_{2}\in
H_{0,L}^{1}(\C)$,
\[
\int_{\C}(v_{2}\Delta v_{1}-v_{1}\Delta
v_{2})\,dxdy=\int_{\Omega}(v_{2}\frac{\partial v_{1}}{\partial
\nu}-v_{1}\frac{\partial v_{2}}{\partial \nu})\,dx,
\]
we see
\[
\int_{\Omega}B_{1/2}g_{2}\cdot
g_{1}\,dx=\int_{\Omega}B_{1/2}g_{1}\cdot g_{2}\,dx
\]
and
\[
\int_{\Omega}v_{2}(x,0)A_{1/2}v_{1}(x,0)\,dx=\int_{\Omega}v_{1}(x,0)A_{1/2}v_{2}(x,0)\,dx.
\]
On the other hand, by using (\ref{eqn-sol-w}) with $\xi=v$ and Lemma
\ref{prop-embedding}, we obtain that $B_{1/2}$ is a positive compact
operator in $L^{2}(\Omega)$. Hence  by the spectral theory for
 self-adjoint compact operators, we have that all the eigenvalues of
$B_{1/2}$ are real, positive, and that  there are corresponding
eigenfunctions which make up an orthonormal basis of
$L^{2}(\Omega)$. Furthermore, such basis and eigenvalues are
explicit in terms of those of the Laplacian with Dirichlet boundary
conditions, since  $(-\Delta)^{-1}$ has $B_{1/2}$ as unique,
positive and self-adjoint square root, by Proposition \ref{prop-B-composition}.
Summarizing:

\begin{prop} \label{lem-linear-eigenvalue} Let $\{\varphi_{k}\}$ be an orthonormal basis of
$L^{2}(\Omega)$ forming a spectral decomposition of $-\Delta$ in
$\Omega$ with Dirichlet boundary conditions, as in
{\rm(\ref{eqn-Dirichlet-egv})}, with $\{\lambda_{k}\}$ the corresponding Dirichlet eigenvalues of
$-\Delta$ in $\Omega$. Then, for all $k\ge 1$,
\begin{equation}\label{eqn-Ahalf-eig}
\left\{
\begin{array}{ll}
A_{1/2}\varphi_{k}=\lambda_{k}^{1/2}\varphi_{k}&\mbox{in}\;\Omega,\\
\varphi_{k}=0 &\mbox{on}\;\partial\Omega.
\end{array}
\right.
\end{equation}
In particular, $\{\varphi_{k}\}$ is also a basis formed by the  eigenfunctions
of $A_{1/2}$, with eigenvalues $\{\lambda_{k}^{1/2}\}$.
\end{prop}

\noindent{\bf Proof of Proposition \ref{prop-V-def}} It follows from
Lemma \ref{lem-V-zero-withd} and Lemma~\ref{lem-spectrum}.

 \hfill{$\Box$}

\medskip

\noindent{\bf Proof of Proposition \ref{prop-A-half-def}} It follows
from Lemma \ref{lem-extension1}, Lemma \ref{lem-spectrum} and its
proof, and   Propositions \ref{prop-B-composition} and
\ref{lem-linear-eigenvalue}.
 \hfill{$\Box$}

\setcounter{equation}{0}
\section{Regularity of solutions}\label{sec-regularity}

In this section we study the regularity of weak solutions for linear
and nonlinear problems involving $A_{1/2}$. First we consider the
linear problem
\begin{equation}\label{eqn-frac-linear}
 \left\{
\begin{array}{ll} A_{1/2}u=g(x) &\mbox{in} \;\Omega,\\
u=0  &\mbox{on}\;\partial\Omega,
\end{array}
\right.
\end{equation}
where $g\in \V_{0}^{*}(\Omega)$ and $\Omega$ is a smooth bounded domain
in $\R^{n}$. By the construction of the previous section, the
precise meaning of (\ref{eqn-frac-linear}) is that
$u=\text{tr}_{\Omega}v$, where the function $v\in H_{0,L}^{1}(\C)$
with $v(\cdot,0)=u\in \V_{0}(\Omega)$ satisfies
\begin{equation}\label{eqn-frac-linear2}
\left\{
\begin{array}{ll}
 \Delta v=0 &  \mbox{in}\;\C, \\
 v=0 & \mbox{on}\;\partial_{L}{\mathcal C},\\
\frac{\partial v}{\partial \nu}=g(x) & \mbox{on} \; \Omega\times
\{0\}.
\end{array}
\right.
\end{equation}
 We will say then that $v$ is a {\it weak
solution} of (\ref{eqn-frac-linear2}) and that $u$ is a  {\it weak
solution} of (\ref{eqn-frac-linear}).

Most of this section contains the proof of the following
analogues of the $W^{2,p}$-estimates of Calder{\'o}n-Zygmund and of the
Schauder estimates.

\begin{prop}\label{prop-linear-reg}
Let $\alpha\in (0,1)$, $\Omega$ be a $C^{2,\alpha}$ bounded domain
of $\R^{n}$, $g\in \V_{0}^{*}(\Omega)$, $v \in H_{0,L}^{1}(\C)$ be the
weak solution of {\rm(\ref{eqn-frac-linear2})}, and
$u=\text{tr}_{\Omega}v$ be the weak solution of {\rm
(\ref{eqn-frac-linear})}. Then,

{\rm(i)} If $g\in L^{2}(\Omega)$, then $u\in
H^{1}_{0}(\Omega)$.

{\rm(ii)} If $g\in H^{1}_{0}(\Omega)$, then $u\in
H^{2}(\Omega)\cap H^{1}_{0}(\Omega)$.

{\rm(iii)} If $g\in L^{\infty}(\Omega)$, then $v\in
W^{1,q}(\Omega\times(0,R))$ for all $R>0$ and $1<q<\infty$. In
particular,
$v\in C^{\alpha}(\overline{\C})$ and  $u\in C^{\alpha}(\overline{\Omega})$.

{\rm(iv)} If  $g\in C^{\alpha}(\overline{\Omega})$ and $g|_{\partial
\Omega}\equiv 0$, then $v\in
  C^{1,\alpha}(\overline{\C})$ and $u\in
  C^{1,\alpha}(\overline{\Omega})$.

   {\rm(v)} If  $g\in C^{1,\alpha}(\overline{\Omega})$ and $g|_{\partial \Omega}\equiv 0$, then $v\in
  C^{2,\alpha}(\overline{\C})$ and $u\in
  C^{2,\alpha}(\overline{\Omega})$.
   \end{prop}

 As a consequence, we deduce the regularity of bounded
 weak solutions to the nonlinear problem
\begin{equation}\label{eqn-frac-non-reg}
 \left\{
\begin{array}{ll} A_{1/2}u=f(u) &\mbox{in} \; \Omega,\\
u=0  &\mbox{on}\;\partial\Omega.
\end{array}
\right.
\end{equation}
As before,
the precise meaning for (\ref{eqn-frac-non-reg}) is that $v\in
H_{0,L}^{1}(\C)$, $v(\cdot,0)=u$, and $v$ is a weak solution of
\begin{equation}\label{eqn-frac-non-reg2}
\left\{
\begin{array}{ll}
 \Delta v=0 &  \mbox{in}\; \C, \\
 v=0 & \mbox{on}\; \partial_{L}{\mathcal C},\\
\frac{\partial v}{\partial \nu}=f(v(\cdot,0))& \mbox{on} \;
\Omega\times \{0\}.
\end{array}
\right.
\end{equation}
Here the weak solution $u$ is assumed to be
bounded. Regularity results for weak solutions not assumed
a priori to be
bounded, of subcritical and critical problems will be proved in
section \ref{sec-h-sub}.

By $C_{0}(\overline{\Omega})$ we denote the space of continuous functions in $\overline{\Omega}$
vanishing on the boundary $\partial \Omega$. In the following result note that
$f(0)=0$ is required to have $C^{1}(\overline{\Omega})$ regularity of solutions
of (\ref{eqn-frac-non-reg}).
\begin{prop}\label{prop-non-reg}
 Let $\alpha\in (0,1)$, $\Omega$ be a $C^{2,\alpha}$ bounded domain of
 $\R^{n}$, and
 $f$ be a $C^{1,\alpha}$ function such that $f(0)=0$.
  If  $u\in L^{\infty}(\Omega)$ is a
 weak solution of {\rm(\ref{eqn-frac-non-reg})},  and thus $v \in H_{0,L}^{1}(\C)\cap
 L^{\infty}(\C)$ is a weak solution of {\rm(\ref{eqn-frac-non-reg2})},
  then  $u\in C^{2,\alpha}(\overline{\Omega})\cap C_{0}(\overline{\Omega})$.
 In addition, $v\in C^{2,\alpha}(\overline{\C})$.
\end{prop}

\noindent{\bf Proof.} By (iii) of
Proposition \ref{prop-linear-reg} we have that $u\in
C^{\alpha}(\overline{\Omega})$. Next, by (iv) of Proposition
\ref{prop-linear-reg} and since on $\partial \Omega\times\{0\}$,
$g:=f(v(\cdot,0))=f(0)=0$, we have $u\in
C^{1,\alpha}(\overline{\Omega})$. Finally, $v\in
C^{2,\alpha}(\overline{\C})$ and $u\in
C^{2,\alpha}(\overline{\Omega})$ from (v) of Proposition
\ref{prop-linear-reg} since $g=f(u)$ vanishes on $\partial\Omega$ and it
is of class $C^{1,\alpha}$,
since both $f$ and $u$ are $C^{1,\alpha}$.

 \hfill{$\Box$}

\medskip

 \noindent{\bf Proof of Proposition \ref{prop-linear-reg}.}
 (i) and (ii).  Both statements follow immediately from Propositions  \ref{prop-V-def}
 and \ref{prop-A-half-def}. Simply use that  $\{\varphi_{k}\}$ is an orthonormal basis of $L^{2}(\Omega)$
 and that
$\{\varphi_{k}/\lambda_{k}^{1/2}\}$ is an orthonormal basis of
$H_{0}^{1}(\Omega)$.  For part (ii), note that if  $A_{1/2}u=g\in
H_{0}^{1}(\Omega)$,  then we have $\Delta u\in L^{2}(\Omega)$.

 (iii) Let $v$ be a weak solution of (\ref{eqn-frac-linear2}).
 We proceed with a useful method, introduced by Cabr{\'e} and Sol\`{a} Morales
  in \cite{CS05}, which consists of  using the auxiliary function
\begin{equation}\label{eqn-f-au}
w(x,y)=\int_{0}^{y}v(x,t)\,dt\;\;\;\text{for}\; (x,y)\in \C.
\end{equation}
 Since $(\Delta w)_{y}=0$ in $\C$, we have that
$\Delta w$ is independent of $y$. Hence we can compute it on
$\{y=0\}$. On $\{y=0\}$, since $w\equiv 0$, we have $\Delta
w=w_{yy}=v_{y}$. Thus $w$ is a solution of the Dirichlet problem
\begin{equation}\label{eqn-l-nonref}
\left\{
 \begin{array}{ll}
-\Delta w(x,y)=g(x) &\mbox{in}\;\C,\\
w=0  &\mbox{on}\; \partial\C.
 \end{array}
 \right.
  \end{equation}

  We extend $w$ to the whole cylinder $\Omega\times\R$ by
  {\it odd} reflection:
  \begin{equation*}
  w_{od}(x,y)=\left\{
\begin{array}{ll}
w(x,y)&\mbox{for}\;y\ge 0,
\\
-w(x,-y)&\mbox{for}\;y\le 0.
\end{array}
  \right.
  \end{equation*}
  Moreover, we put
\begin{equation*}
 g_{od}(x,y)=\left\{
\begin{array}{ll}
g(x)&\mbox{for}\;y> 0,
\\
-g(x)&\mbox{for}\;y< 0.
\end{array}
  \right.
  \end{equation*}
Then we obtain
\begin{equation}\label{eqn-reflection-y}
\left\{
 \begin{array}{ll}
-\Delta w_{od}=g_{od}&
 \mbox{in}\; \Omega\times\R,\\
w_{od}=0& \mbox{on}\;\partial \Omega\times\R.
 \end{array}
 \right.
  \end{equation}
  Since $g_{od}\in L^{q}(\Omega\times(-2R,2R))$ for all $R>0$ and $1<q<\infty$,
   regularity for the Dirichlet problem (\ref{eqn-reflection-y}) gives $w_{od}\in
  W^{2,q}(\Omega\times(-R,R))$
   for all $R>0$ and $1<q<\infty$. In particular,
  $w\in C^{1,\alpha}(\overline{\C})$.
Therefore, $v=w_{y}\in C^{\alpha}(\overline{\C})$ and $u\in C^{\alpha}(\overline{\Omega})$.

(iv) Choose a smooth domain $H$ such that $\overline{\Omega}\subset
H$, and let
\[
g_{H}=\left\{
\begin{array}{ll}
g &\text{in}\; \overline{\Omega},
\\
0 & \text{in}\; H \setminus\overline{\Omega}.
\end{array}
\right.
\]
 We have that $g_{H}\in C^{\alpha}(\overline{H})$,
 since $g\!\!\mid_{\partial \Omega}=0$, by assumption.
 Consider the weak solution $v_{H}$ of
\begin{equation*}
\left\{
\begin{array}{ll}
 \Delta v_{H}=0 &  \mbox{in}\quad H\times (0,\infty), \\
 v_{H}=0 & \mbox{on}\quad \partial H\times [0,\infty),\\
\frac{\partial v_{H}}{\partial \nu}=g_{H}(x) & \mbox{on} \quad
H\times \{0\}.
\end{array}
\right.
\end{equation*}
Consider also the auxiliary function
 \[
w_{H}(x,y)=\int_{0}^{y}v_{H}(x,t)\,dt\;\;\;\text{in}\;
\overline{H}\times[0,\infty),
\]
which solves problem (\ref{eqn-l-nonref}) with $\Omega$ and $g$
replaced by $H$ and $g_{H}$.

 Using boundary regularity theory (but away from the corners of $H\times[0,\infty)$)
 for this Dirichlet problem, we see
 that  $w_{H}$ is $C^{2, \alpha}(H\times (0,\infty))$ (again, here we
 do not claim regularity at the corners
 $\partial H\times \{0\}$). Thus, $w_{H}\in
C^{2, \alpha}(\overline{\C})$ (here instead we include the corners
$\partial \Omega\times\{0\}$ of $\C$).

Consider the difference $\varphi=w_{H}-w$ in $\C$, where $w$ is
defined by (\ref{eqn-f-au}). It is clear that
\begin{equation*}
\left\{
\begin{array}{ll}
 \Delta \varphi=0 &  \mbox{in}\quad \C, \\
 \varphi=w_{H} & \mbox{on}\quad \partial_{L}{\mathcal C},\\
\varphi=0& \mbox{on} \quad \Omega\times \{0\}.
\end{array}
\right.
\end{equation*}
 We extend $\varphi$ to the whole cylinder $\Omega\times\R$ by
  {\it odd} reflection:
  \begin{equation*}
  \varphi_{od}(x,y)=\left\{
\begin{array}{ll}
\varphi(x,y)&\mbox{for}\;y\ge 0,
\\
-\varphi(x,-y)&\mbox{for}\;y\le 0.
\end{array}
  \right.
  \end{equation*}
  Moreover, we put
\begin{equation*}
  w_{H, od}(x,y)=\left\{
\begin{array}{ll}
w_{H}(x,y)&\mbox{for}\;y> 0,
\\
-w_{H}(x,-y)&\mbox{for}\;y\le 0.
\end{array}
  \right.
  \end{equation*}
Then we have
\begin{equation}\label{eqn-varphi-od}
\left\{
\begin{array}{ll}
 \Delta \varphi_{od}=0 &  \mbox{in}\quad \Omega\times\R, \\
 \varphi_{od}=w_{H, od} & \mbox{on}\quad \partial\Omega\times\R.\\
\end{array}
\right.
\end{equation}

Since $w_{H}\in C^{2,\alpha}(\overline{\C})$, $w_{H}\equiv 0$ on
$\partial \Omega\times\{0\}$, and
$\partial_{yy}w_{H}=\partial_{y}v_{H}=-g_{H}=-g=0$ on $\partial\Omega\times
\{0\}$, we deduce that $w_{H,od}\in
C^{2,\alpha}(\partial\Omega\times \R)$. It follows  from elliptic
regularity for (\ref{eqn-varphi-od}) that $\varphi_{od}\in C^{2,
\alpha}(\overline{\Omega}\times \R)$. Thus, $\varphi\in C^{2,
\alpha}(\overline{\C})$, $w\in C^{2,
\alpha}(\overline{\C})$ and $v=\partial_{y}w\in C^{1, \alpha}(\overline{\C})$.

(v) Choose a smooth bounded domain $B$ such that $\Omega
\subset\overline{B}$. $B$ could be the same as $H$ in (ii), for
instance a ball, but we change its name for notation clarity. Since
$g\in C^{1,\alpha}(\overline{\Omega})$, there exists an extension
$g_{B}\in C^{1,\alpha}(\overline{B})$; see \cite{GT01}. Consider the
solution $v_{B}$ of
\begin{equation*}
\left\{
\begin{array}{ll}
 \Delta v_{B}=0 &  \mbox{in}\quad B\times (0,\infty), \\
 v_{B}=0 & \mbox{on}\quad \partial B\times [0,\infty),\\
\frac{\partial v_{B}}{\partial \nu}=g_{B} & \mbox{on} \quad B\times
\{0\}.
\end{array}
\right.
\end{equation*}
 Consider
the auxiliary function
 \[
w_{B}(x,y)=\int_{0}^{y}v_{B}(x,t)\,dt\;\;\;\text{in}\;
\overline{B}\times[0,\infty).
\]

As before, from interior boundary regularity for the Dirichlet problem of the type
\eqref{eqn-l-nonref} satisfied by $w_{B}$, we obtain that
$w_{B}\in C^{3,\alpha}(B\times [0,\infty))$ since $g_{B}\in
C^{1,\alpha}(\overline{B})$ (away from the corners $\partial
B\times\{0\}$). Thus, $v_{B} \in C^{2,\alpha}(B\times [0,\infty))$.
Thus, $v_{B}\in C^{2,\alpha}(\overline{\C})$. Consider the
difference $\psi=v_{B}-v$ in $\C$, where $v$ is a weak solution of
(\ref{eqn-frac-linear2}). We have that $\psi=v_{B}-v$ satisfies
\begin{equation*} \left\{
\begin{array}{ll}
 \Delta \psi=0 &  \mbox{in}\quad \C, \\
 \psi=v_{B} & \mbox{on}\quad \partial_{L}{\mathcal C},\\
\frac{\partial \psi}{\partial \nu}=0& \mbox{on} \quad \Omega\times
\{0\}.
\end{array}
\right.
\end{equation*}

 We extend $\psi$ to the whole cylinder $\overline{\Omega}\times\R$ now by
  {\it even} reflection:
  \begin{equation*}
  \psi_{ev}(x,y)=\left\{
\begin{array}{ll}
\psi(x,y)&\mbox{for}\;y\ge 0,
\\
\psi(x,-y)&\mbox{for}\;y\le 0.
\end{array}
  \right.
  \end{equation*}
  Moreover, we put
\begin{equation*}
  v_{B,ev}(x,y)=\left\{
\begin{array}{ll}
v_{B}(x,y)&\mbox{for}\;y> 0,
\\
v_{B}(x,-y)&\mbox{for}\;y\le 0.
\end{array}
  \right.
  \end{equation*}
Then, since $\frac{\partial \psi}{\partial \nu}=0$ on
$\Omega\times\{0\}$, we have
\begin{equation*}
\left\{
\begin{array}{ll}
 \Delta \psi_{ev}=0 &  \mbox{in}\quad \Omega\times\R, \\
 \psi_{ev}=v_{B, ev} & \mbox{on}\quad \partial\Omega\times\R.\\
\end{array}
\right.
\end{equation*}
Since $v_{B}\in C^{2,\alpha}(\overline{\C})$,
$-\partial_{y}v_{B}=g_{B}=g=0$ on $\partial\Omega\times \{0\}$, we deduce that $v_{B,ev}\in  C^{2,
\alpha}(\partial\Omega\times \R)$. Therefore, it follows from classical
regularity that $\psi_{ev}\in C^{2, \alpha}(\overline{\Omega}\times
\R)$. Thus, $\psi\in C^{2, \alpha}(\overline{\C})$, and $v\in C^{2,
\alpha}(\overline{\C})$. \hfill{$\Box$}

\setcounter{equation}{0}
\section{Maximum principles}\label{sec-h-max}
In this section we establish several maximum principles for
$A_{1/2}$. We denote by $C_{0}(\overline{\Omega})$ the space of continuous
functions in $\overline{\Omega}$ vanishing on the boundary $\partial \Omega$. For
convenience, we state the results for functions in
$C_{0}(\overline{\Omega})\cap C^{2}(\overline{\Omega})$ (a space
contained in  $H_{0}^{1}(\Omega)\subset \V_{0}(\Omega)$), but this
can be weakened.

The first statement is the weak maximum principle.
\begin{lema}\label{lem-max1}
Assume that $u\in C^{2}(\overline{\Omega})$ satisfies
\begin{equation*}
\left\{
\begin{array}{ll}
A_{1/2}u +c(x)u \ge 0 &\mbox{in}\; \Omega,\\
u=0 &\mbox{on}\; \partial \Omega,
\end{array}
\right.
\end{equation*}
where $\Omega$ is a smooth bounded domain in $\R^{n}$ and $c\ge 0$
in $\Omega$. Then, $u\ge 0$ in $\Omega$.
\end{lema}
\noindent{\bf Proof.} Consider the extension $v=\mbox{h-ext}(u)$. If
we prove that $v\ge 0$ in $\C$, then $u\ge 0$ in $\Omega$. Suppose
by contradiction that $v$ is negative somewhere in $\C$. Then, since
$\Delta v=0$ in $\C$ and $v=0$ on $\partial_{L}\C$, we deduce that
$v$ is negative somewhere in $\Omega\times\{0\}$ and that $\inf_{\C}
v <0$ is achieved at some point $(x_{0},0)\in \Omega\times\{0\}$.
Thus, we have
\[
\inf_{\C}v=v(x_{0},0)<0.
\]
By Hopf's lemma,
\[
v_{y}(x_{0},0) >0.
\]
It follows
\[
\frac{\partial v}{\partial \nu}=-v_{y}(x_{0},0)=A_{1/2}v(x_{0},0)
<0.
\]
Therefore, since $c\ge 0$,
\[
A_{1/2}v(x_{0},0)+c(x_{0})v(x_{0},0)<0.
\]
 This is a contradiction  with the hypothesis $A_{1/2} u+ c(x)u\ge 0$.
 \hfill {$\Box$}

\medskip
The next statement is the strong maximum principle for $A_{1/2}$.
\begin{lema}  \label{lem-max2}
Assume that $u\in
C^{2}(\overline{\Omega})$ satisfies
\begin{equation*}
\left\{
\begin{array}{ll}
A_{1/2}u +c(x)u \ge 0 &\mbox{in}\; \Omega,\\
u\ge 0 &\mbox{in}\;  \Omega,\\
u=0 &\mbox{on}\; \partial \Omega,
\end{array}
\right.
\end{equation*}
where $\Omega$ is a smooth bounded domain in $\R^{n}$ and $c\in
L^{\infty}(\Omega)$. Then, either $u> 0$ in $\Omega$, or $u\equiv 0$
in $\Omega$.
\end{lema}
\noindent{\bf Proof.} The proof is similar to that of Lemma
\ref{lem-max1}. Consider $v=\mbox{h-ext}(u)$. We observe that $v\ge
0$ in $\C$. Suppose that $v\not\equiv 0$ but $u=0$ somewhere in
$\Omega$. Then there exists a minimum point
$(x_{0},0)\in\Omega\times\{0\}$ of $v$ where $v(x_{0},0)=0$. Then by
Hopf's lemma we see that $ A_{1/2}u(x_{0})=- v_{y}(x_{0},0)<0$. This
implies that $A_{1/2}u(x_{0})+c(x_{0})u(x_{0})<0$, because of
$v(x_{0},0)=u(x_{0})=0$.
 \hfill {$\Box$}

\medskip

Next we establish a Hopf lemma for $A_{1/2}$, following a proof from \cite{CCFS}.
 \begin{lema}  \label{lem-hopf}
Let $\Omega$ be a bounded domain in $\R^{n}$ and $c\in
L^{\infty}(\Omega)$.

{\rm (i)} Assume that $\Omega$ is smooth and that  $
0\not\equiv u\in
C^{2}(\overline{\Omega})$ satisfies
\begin{equation*}
\left\{
\begin{array}{ll}
A_{1/2}u +c(x)u \ge 0 &\mbox{in}\; \Omega,\\
u\ge 0 &\mbox{in}\;  \Omega,\\
 u=0 &\mbox{on}\; \partial \Omega.
\end{array}
\right.
\end{equation*}
Then, $\frac{\partial u}{\partial \nu_{0}}< 0$ on $\partial \Omega$,
where $\nu_{0}$
 is the unit outer normal to $\partial\Omega$.

 {\rm (ii)} Assume that $P\in \partial \Omega$ and that $\partial \Omega$
 is smooth in a neighborhood of~$P$. Let $0\not\equiv v \in C^{2}(\overline{\C})\cap
L^{\infty}(\C)$, where $\C=\Omega\times(0,\infty)$, satisfy
\begin{equation*}
\left\{
\begin{array}{ll}
\Delta v=0  &\mbox{in}\; \C,\\
v\ge 0 &\mbox{on}\; \partial_{L} \C,\\
\frac{\partial v}{\partial \nu}+c(x)v\ge 0 &\mbox{on}\;
\Omega\times\{0\}.
\end{array}
\right.
\end{equation*}
If $v(P,0)=0$, then $\frac{\partial v(P,0)}{\partial
\nu_{0}}<0$, where $\nu_{0}$
 is the unit outer normal in~$\R^{n}$ to $\partial\Omega$.
\end{lema}

\noindent{\bf Proof.} We follow the  proof given in \cite{CCFS}.
Note that
statement (i) is a particular case of (ii). Thus, we only need to prove (ii).

{\it Step 1.} We shall first prove the lemma in the case $c\equiv
0$.  Without loss of generality we may assume that
$(P,0)=P_{1}=(b_{1},0,\cdots,0)\in \partial \Omega \times\{0\}$,
$b_{1}>0$ and $\nu_{0}=(1,0,\cdots,0)$. Hence we need to prove
\[
\frac{\partial v(P_{1})}{\partial x_{1}}<0.
\]
Since $\Omega$ is smooth in a neighborhood of $P$, there is a
half-ball in $\R^{n+1}_{+}$ included in the domain $\C$, such that
$P_{1}$ is the only point in the closed half-ball belonging also to
$\partial_{L} \C$. Let  $P_{2}\in \Omega\times\{0\}$ and $r>0$ be the center
and radius of such ball. Then we have $P_{2}=(b_{2},0,\cdots,0)\in
\Omega\times\{0\}$. Denote
\begin{equation*}
\begin{split}
& B_{r}^{+}(P_{2}):=\{z=(x,y)\mid |z-P_{2}|< |P_{1}-P_{2}|=:r, y>0  \}\subset \C, \\
& B_{r/2}^{+}(P_{2}):=\{z=(x,y)\mid |z-P_{2}|< |P_{1}-P_{2}|/2, y>0  \},\\
& A=B_{r}^{+}(P_{2})\setminus \overline{B_{r/2}^{+}(P_{2})}.
\end{split}
\end{equation*}
Recall that $P_{1}\in \partial
B_{r}^{+}(P_{2})\cap(\partial\Omega\times\{0\})$.

Consider the function on $A$:
\[
\varphi(z)=\exp (-\lambda |z-P_{2}|^{2})-\exp (-\lambda
|P_{1}-P_{2}|^{2}),
\]
with $\lambda >0$ to be determined later. Note that
\[
\Delta \varphi=\exp (-\lambda
|z-P_{2}|^{2}) \left\{ 4\lambda^{2}|z-P_{2}|^{2}-2(n+1)\lambda\right\}.
\]
We can choose $\lambda>0$ large enough  such that $\Delta \varphi\ge
0$ in $A$.

On the other hand, by Lemma \ref{lem-max2}, we see that $v>0$ in
$\overline{A}\setminus \{P_{1}\}$.  Hence, since $\varphi\equiv 0$
on $\partial B_{r}^{+}(P_{2})\cap \{y>0\}$, we can take
$\varepsilon>0$ such that
\[
v-\varepsilon \varphi \ge 0 \quad \mbox{on} \; \partial A\cap
\{y>0\}.
\]
Since $-\Delta (v-\varepsilon \varphi)\ge 0$ in $A$, and
\[-\partial_{y}(v-\varepsilon\varphi)=\frac{\partial v}{\partial
\nu}\ge 0\quad\text{on}\; \partial A\cap \{y=0\},\] (recall that
$c\equiv 0$)
 by the maximum
principle as in Lemma \ref{lem-max1}
 we obtain
\[
v-\varepsilon \varphi \ge 0 \quad \mbox{in} \;  A.
\]
Thus, from $v-\varepsilon \varphi=0$ at $P_{1}$ we see that
$\partial_{x_{1}} (v-\varepsilon \varphi)(P_{1})\le
0$. Therefore, $\partial_{ x_{1}} v(P_{1})\le
\varepsilon\partial_{x_{1}} \varphi(P_{1})=-2\lambda(b_{1}-b_{2})e^{-\lambda|P_{1}-P_{2}|^{2}}< 0$. Thus
 we have the desired result.

{\it Step 2}. In the case $c\not\equiv 0$, we define the function
$w=v\exp (-\beta y)$ for some $\beta>0$ to be determined. From a
direct calculation, we see that
\[
-\Delta w-2\beta\partial_{y} w= \beta^{2}w\ge 0 \quad \mbox{in} \;
\C
\]
and, choosing $\beta \ge \|c\|_{L^{\infty}(\Omega)}$,
\[
-\partial_{y}w\ge [\beta-c(x)]w\ge 0 \quad \mbox{on} \;
\Omega\times\{0\}.
\]
 Now we can apply to $w$ the same approach  as in Step 1, with $\Delta$
replaced by $\Delta+2\beta\partial_{y}$, and obtain the assertion.
\hfill {$\Box$}

 \medskip

 Finally, we establish a maximum principle for $A_{1/2}$ in domains of small measure. Note that
 in part (ii) of its statement, the hypothesis on  small measure  is made only on the base of $\Omega$
 of the cylinder $\C$.

 \begin{prop}  \label{lem-max3}
{\rm(i)} Assume that $u\in
C^{2}(\overline{\Omega})$ satisfies
\begin{equation*}
\left\{
\begin{array}{ll}
A_{1/2}u +c(x)u \ge 0 &\mbox{in}\; \Omega,\\
u=0 &\mbox{on}\; \partial \Omega,
\end{array}
\right.
\end{equation*}
where $\Omega$ is a smooth bounded domain in $\R^{n}$ and  $c\in
L^{\infty}(\Omega)$. Then, there exists $\delta>0$ depending only on
$n$ and $\|c^{-}\|_{L^{\infty}(\Omega)}$, such that if
$|\Omega\cap\{u<0\}|\le \delta$, then $u\ge 0$ in $\Omega$.

{\rm(ii)} Assume that $\Omega$ is a bounded (not necessary smooth)
domain of $\R^{n}$ and $c\in L^{\infty}(\Omega)$. Let $v\in
C^{2}(\overline{\C})\cap L^{\infty}(\C)$, where
$\C=\Omega\times(0,\infty)$, satisfy
\begin{equation*}
\left\{
\begin{array}{ll}
\Delta v=0  &\mbox{in}\; \C,\\
v\ge 0 &\mbox{on}\; \partial_{L} \C,\\
\frac{\partial v}{\partial \nu}+c(x)v\ge 0 &\mbox{on}\;
\Omega\times\{0\}.
\end{array}
\right.
\end{equation*}
Then, there exists $\delta>0$ depending only on $n$ and
$\|c^{-}\|_{L^{\infty}(\Omega)}$, such that if  $|\Omega\cap
\{v(\cdot,0)<0\}|\le \delta$ then $v\ge 0$ in $\C$.
\end{prop}
\noindent{\bf Proof.} For part (i) of the theorem, consider
$v=\mbox{h-ext}(u)$. We see that $v$ satisfies the assumptions on
part (ii) of the theorem. Hence, it is enough to prove part (ii).
For this, let $v^{-}=\max\{0,-v\}\ge 0$. Since $v^{-}=0$ on
$\partial \Omega\times [0,\infty)$, we see
\[
0=\int_{\C}v^{-}\Delta
v\,dxdy=\int_{\Omega\times\{0\}}v^{-}\frac{\partial
v}{\partial\nu}\,dx+ \int_{\C}|\nabla v^{-}|^{2}\,dxdy.
\]
Then,
\begin{equation*}
\begin{split}
\int_{\C}|\nabla v^{-}|^{2}\,dxdy=&-\int_{\Omega\times\{0\}}v^{-}\frac{\partial v}{\partial\nu}\,dx\\
\le & \int_{\Omega\times\{0\}}v^{-}cv\,dx=
\int_{\Omega}-c(v^{-})^{2}\,dx
\\ \le&  \int_{\Omega\cap\{v^{-}(\cdot,0)>0\}}c^{-}(v^{-}(\cdot,0))^{2}\,dx \\
\le &
|\Omega\cap\{v^{-}(\cdot,0)>0\}|^{1/n}\|c^{-}\|_{L^{\infty}(\Omega)}\|v^{-}(\cdot,0)\|_{L^{2n/(n-1)}(\Omega)}^{2}.
\end{split}
\end{equation*}
Thus,  extending $v^{-}$ by $0$ outside $\overline{\C}$ we obtain an
$H^{1}(\R^{n+1}_{+})$ function and thus we have
\begin{equation*}
\begin{split}
0<S_{0}\le&\frac{\int_{\R^{n+1}_{+}}|\nabla
v^{-}|^{2}\,dxdy}{\|v^{-}(\cdot,0)\|_{L^{2n/(n-1)}(\R^{n})}^{2}} =
\frac{\int_{\C}|\nabla
v^{-}|^{2}\,dxdy}{\|v^{-}(\cdot,0)\|_{L^{2n/(n-1)}(\Omega)}^{2}}\\
\le &
|\Omega\cap\{v^{-}(\cdot,0)>0\}|^{1/n}\|c^{-}\|_{L^{\infty}(\Omega)},
\end{split}
\end{equation*}
 where $S_{0}$ is the best constant of the
Sobolev trace inequality in $\R^{n+1}_{+}$. If
$|\Omega\cap\{v^{-}(\cdot,0)>0\}|$ is small enough, we arrive at a
contradiction.
 \hfill {$\Box$}


\setcounter{equation}{0}
\section[The subcritical case]{Subcritical case
and $L^{\infty}$ estimate of Brezis-Kato type}\label{sec-h-sub} In
this section, we study the nonlinear problem (\ref{eqn-frac}) with
$f(u)=u^{p}$ in the subcritical and critical cases. In the subcritical case we look for a function
$v(x,y)$ satisfying for $x\in \Omega$ and $y\in\R^{+}$,
\begin{equation}\label{eqn-B}
\left\{
\begin{array}{ll}
 \Delta v=0 &  \mbox{in}\quad {\mathcal C}=\Omega\times (0,\infty),\\
 v=0 & \mbox{on}\quad \partial_{L}{\mathcal C}=\partial\Omega\times [0,\infty),\\
\frac{\partial v}{\partial \nu}=v^{p} & \mbox{on} \quad \Omega\times\{0\},\\
v>0 &  \mbox{in}\quad {\mathcal C},
\end{array}
\right.
\end{equation}
where $\nu$ is the unit outer normal to $\C$ at $\Omega\times\{0\}$
and $1< p<2^{\sharp}-1$ if $n\ge 2$, or $1< p<\infty$ if $n=1$. If
$v$ is a solution of \eqref{eqn-B}, then $v(x,0)=u(x)$ is a solution
of \eqref{eqn-frac} with the nonlinearity $f(u)=u^{p}$.

In order to find a solution of \eqref{eqn-B} as stated in Theorem
\ref{teo-sub1}, we consider the following minimization problem:
\[
I_{0}=\inf\left\{ \int_{{\mathcal{C}}}|\nabla v(x,y)|^{2}dxdy \mid
v\in H_{0,L}^{1}(\C), \int_{\Omega}|v(x,0)|^{p+1}dx=1\, \right\}.
\]
We show that $I_{0}$ is achieved.
\begin{prop}\label{prop-achieved1} Assume that $1< p<2^{\sharp}-1$
if $n\ge 2$ or $1< p<\infty$ if $n=1$. Then $I_{0}$ is achieved in
$H_{0,L}^{1}(\C)$ by a nonnegative function $v$.
\end{prop}
\noindent{\bf Proof.} First, there is a function $v\in
H_{0,L}^{1}(\C)$ such that
\[
 \int_{{\mathcal{C}}}|\nabla v(x,y)|^{2}dxdy <\infty\;\; \mbox{and} \;
 \int_{\Omega}|v(x,0)|^{p+1}dx=1.
 \]
 In fact, it suffices to take any $C^{\infty}$ function with compact support in
 $\Omega\times[0,\infty)$ and not identically zero on $\Omega\times\{0\}$, and multiply it by an
 appropriate constant.
  Next we complete the proof by weak lower semi-continuity of the Dirichlet integral
  and by the compact embedding property in Lemma \ref{prop-embedding}.
  Finally, note that
  $|v|\ge 0$ is a nonnegative minimizer if $v$ is a minimizer.
 \hfill {$\Box$}

 \medskip

To establish the regularity of the minimizer just obtained,
we prove an $L^{\infty}$-estimate of
Brezis-Kato type by the technique of bootstrap for subcritical or
critical nonlinear problems. Let
 $g_{0}$ be a Carath\'{e}odory function in $\Omega\times \R$ satisfying the growth
condition
\begin{equation}\label{eqn-h-gzero}
|g_{0}(x, s)|\le C(1+|s|^{p})\quad\text{for all}\; (x,s)\in
\Omega\times\R,
\end{equation}
where $\Omega$ is a smooth domain in $\R^{n}$, $1\le p\le
\frac{n+1}{n-1}$ if $n\ge 2$, or  $1\le p<\infty$ if $n=1$. We
consider the problem
\begin{equation}\label{eqn-boot}
\left\{
\begin{array}{ll}
 \Delta v=0 &  \mbox{in}\quad {\mathcal C}=\Omega\times (0,\infty),\\
 v=0 & \mbox{on}\quad \partial_{L}{\mathcal C}=\partial\Omega\times [0,\infty),\\
\frac{\partial v}{\partial \nu}=g_{0}(\cdot, v)& \mbox{on} \quad
\Omega\times \{0\}.
\end{array}
\right.
\end{equation}

\begin{teo} \label{prop-teo-bk}
 Let $v\in H_{0,L}^{1}(\C)$ be a weak solution of
 {\rm(\ref{eqn-boot})}
and assume the growth condition {\rm(\ref{eqn-h-gzero})} for
$g_{0}$, with $1\le p\le \frac{n+1}{n-1}$ if $n\ge 2$, or  $1\le
p<\infty$ if $n=1$.
  Then,   $v(\cdot,0)\in L^{\infty}(\Omega)$.
\end{teo}
\noindent {\bf Proof.} The proof follows the one of Brezis-Kato for
the Laplacian. First of all, let us rewrite the condition on $g_{0}$
as
\[
|g_{0}(x, v)|\le a(x)(1+|v(x,0)|)
\]
with a function
\[
a(x):=\frac{|g_{0}(x,v(x,0))|}{1+|v(x,0)|}
\]
which satisfies
\[
0\le a\le C(1+|v(\cdot,0)|^{p-1})\in L^{n}(\Omega),
\]
since $v\in H_{0,L}^{1}(\C)$, $v(\cdot,0)\in
L^{\frac{2n}{n-1}}(\Omega)$ and $p-1\le \frac{2}{n-1}$.

 Denote
 \[
B_{r}^{+}=\{(x,y)\mid |(x,y)|< r\;\text{and}\; y>0 \}.
\]
   For $\beta \ge 0$ and
$T>1$, let $\varphi_{\beta,\,T}=v v_{T}^{2\beta}\in H_{0,L}^{1}(\C)$ with $v_{T}=\min\{|v|,T\}$. Denote
\[
  D_{T}=\{(x,y)\in\C\mid|v(x,y)|< T\}.
\]
By direct computation, we see
\begin{align*}
\int_{\C}|\nabla
(vv_{T}^{\beta})|^{2}\,dxdy=&\int_{\C}v_{T}^{2\beta}|\nabla
v|^{2}\,dxdy\\&+\int_{D_{T}}(2\beta+\beta^{2})|v|^{2\beta}|\nabla
v|^{2}\,dxdy.
\end{align*}
Multiplying (\ref{eqn-boot}) by $\varphi_{\beta,\,T}$ and
integrating by parts, we obtain
\begin{align*}
&\int_{\C}v_{T}^{2\beta}|\nabla v|^{2}\,dxdy+2\beta\int_{D_{T}}|v|^{2\beta}|\nabla v|^{2}\,dxdy\\
= &\int_{\C}\nabla v\nabla ( vv_{T}^{2\beta})\,dxdy
=\int_{\Omega\times\{0\}}g_{0}(x,v)vv_{T}^{2\beta}\,dx\\
\le& \int_{\Omega\times\{0\}}a(x)(1+|v|)^{2}v_{T}^{2\beta}\,dx.
\end{align*}
Combining these facts, we have
\begin{align*}
\int_{\C} |\nabla (vv_{T}^{\beta})|^{2} \,dxdy \le&
C(\beta+1)\int_{\Omega\times\{0\}}a(x)(1+|v|^{2})v_{T}^{2\beta}\,dx,
\end{align*}
where $C$ denotes different constants independent of $T$ and of
$\beta$.
 By Lemma \ref{p-in}, we deduce
\begin{equation}
\begin{split}
\label{boot-s-le-two} \Big(\int_{\Omega\times\{0\}}
|vv_{T}^{\beta}|^{2^{\sharp}} \,dx \Big)^{2/2^{\sharp}}\le&
C(\beta+1)\int_{\Omega\times\{0\}}a(x)(1+|v|^{2})v_{T}^{2\beta}\,dx.
\end{split}
\end{equation}

Assume that $|v(\cdot,0)|^{\beta+1}\in L^{2}(\Omega)$ for some
$\beta\ge 0$. Then we obtain that
$\int_{\Omega\times\{0\}}|v|^{2}v_{T}^{2\beta}\,dx$ and
$\int_{\Omega\times\{0\}}v_{T}^{2\beta}\, dx$ are  bounded uniformly
in $T$. In what follows, let $C$ denote constants independent of $T$
---but that may depend on $\beta$ and
$\|v(\cdot,0)^{\beta+1}\|_{L^{2}(\Omega)}$. Given $M_{0}>0$, we have
\begin{equation*}
\begin{split}
\int_{\Omega\times\{0\}}a|v|^{2}v_{T}^{2\beta}\,dx &\le
M_{0}\int_{\Omega\times\{0\}}|v|^{2}v_{T}^{2\beta}\,dx+\int_{\{a\ge
M_{0}\}}a|v|^{2}v_{T}^{2\beta}\,dx\\
 &\le CM_{0}+\Big(\int_{\{a\ge M_{0}\}}a^{n}\,dx\Big)^{1/n}
\Big(\int_{\Omega\times\{0\}}|vv_{T}^{\beta}|^{2^{\sharp}}\,dx\Big)^{2/2^{\sharp}}\\
&\le  CM_{0}+\varepsilon(M_{0})
\Big(\int_{\Omega\times\{0\}}|vv_{T}^{\beta}|^{2^{\sharp}}\,dx\Big)^{2/2^{\sharp}},
\end{split}
\end{equation*}
where $\varepsilon(M_{0})=(\int_{\{a\ge
M_{0}\}}a^{n}\,dx)^{1/n}\rightarrow 0$ as $M_{0}\rightarrow \infty$.
Note that we can deal with $\int_{\Omega\times\{0\}}a
v_{T}^{2\beta}\,dx$ in the analogue procedure. Therefore, we deduce
from the last inequalities and (\ref{boot-s-le-two}), taking $M_{0}$
large enough so that $C(\beta+1)\varepsilon(M_{0})=\frac{1}{2}$, that
\begin{equation}\label{eqn-iteration}
\Big(\int_{\Omega\times\{0\}}  |vv_{T}^{\beta}|^{2^{\sharp}} \,dx
\Big)^{2/2^{\sharp}}\le C(1+M_{0}).
\end{equation}
 Thus letting $T\rightarrow \infty$, since  $C$ is independent of $T$,
 we obtain that $|v(\cdot,0)|^{\beta+1}\in
L^{2^{\sharp}}(\Omega)$. This conclusion followed simply from
assuming $|v(\cdot,0)|^{\beta+1}\in L^{2}(\Omega)$.

Hence, by iterating $\beta_{0}=0$,
$\beta_{i}+1=(\beta_{i-1}+1)\frac{n}{n-1}$ if $i\ge 1$ in
(\ref{eqn-iteration}), we conclude that  $v(\cdot,0)\in
L^{q}(\Omega)$ for all $q<\infty$. Finally, the proof of part (iii) in
Proposition \ref{prop-linear-reg}  ---which only uses $g\in
L^{q}(\Omega)$ for all $q<\infty$ and not $g\in
L^{\infty}(\Omega)$--- applied with $g(x)=g_{0}(x,v(x,0))$, which
satisfies $|g|\le C(1+|v(\cdot,0)|^{p})\in L^{q}(\Omega)$ for all
$q<\infty$, leads to $v(\cdot,0)\in
C^{\alpha}(\overline{\Omega})\subset L^{\infty}(\Omega)$.
 \hfill {$\Box$}

\medskip

\noindent {\bf Proof of Theorem \ref{teo-sub1}.} Proposition
\ref{prop-achieved1} gives the existence of a weak nonnegative
solution $v$ to (\ref{eqn-B}) after multiplying the nonnegative
minimizer of $I_{0}$ by a constant to take care of the Lagrange
multiplier. Then, Theorem \ref{prop-teo-bk} gives that $v(\cdot,0)\in
L^{\infty}(\Omega)$. Next, Proposition~\ref{prop-non-reg} gives that
$u\in C^{2,\alpha}(\overline{\Omega})$, since $f(s)=|s|^{p}$ is a
$C^{1,\alpha}$ function for some $\alpha \in (0,1)$. Finally, the
strong maximum principle, Lemma \ref{lem-max2}, leads to  $u>0$ in
$\Omega$.
 \hfill{$\Box$}

\setcounter{equation}{0}
\section[A priori estimates]{A priori estimates for positive solutions}\label{sec-h-apriori}
In this section we prove Theorem \ref{teo-ap}. Namely, we establish
 a priori estimates of Gidas-Spruck type for weak solutions of
\begin{equation}\label{eqn-bd1}
\begin{cases}
\Delta v=0 & \mbox{in}\;\C=\Omega\times (0,\infty) \subset \R^{n+1}_{+}, \\
v=0&\mbox{on}\;\partial_{L}\C=\partial \Omega\times[0,\infty),\\
\frac{\partial v}{\partial \nu}=v^{p}& \mbox{on}\; \Omega\times\{0\},\\
v>0&\mbox{in}\; \C,
\end{cases}
\end{equation}
where $\Omega\subset \R^{n}$ is a bounded smooth domain, $n\ge 2$,
and $1< p<\frac{n+1}{n-1}$.

For this, we need  two nonlinear Liouville theorems for problems
involving the square root of the Laplacian in unbounded domain
---one in the whole space, another in the half-space. The
first one was proved by Y.Y. Li, Zhang and Zhu in \cite{LiZhu95}, \cite{LiZ} and Ou in
\cite{Ou}. Its statement is the following ---and it is equivalent to Theorem \ref{teo-LZO-half}
in the Introduction.
\begin{teo}  \label{li-ou} {\rm (\cite{LiZhu95}, \cite{LiZ}, \cite{Ou})}
 For $n\ge 2$ and   $1< p<2^{\sharp}-1=\frac{n+1}{n-1}$,
there exists no weak solution of problem
\begin{equation}\label{eqn-li-ou}
\begin{cases}
\Delta v=0 &\text{in}\; \R^{n+1}_{+},\\
\frac{\partial v }{\partial \nu} =v^{p} & \text{on}\; \partial\R^{n+1}_{+}, \\
v >0  &\text{in}\; \R^{n+1}_{+}.
\end{cases}
\end{equation}
\end{teo}

We need to prove an analogue nonlinear Liouville type result
involving the square root of $-\Delta$ with Dirichlet boundary value
in the half-space. This
is Theorem \ref{thm-l-d2} of the Introduction and
Proposition \ref{prop-sol-d2} in this section. As we will see, this nonlinear Liouville theorem
 in $\R_{+}^{n}$ will be first reduced to the one
dimensional case $\R_{+}$, by using the moving planes method.  After
this,  we prove that there exists no positive bounded solution for
the nonlinear Neumann boundary problem in the quarter $\R^{2}_{++}$,
which corresponds to the nonlinear Liouville theorem involving the
square root of $-\Delta$ with Dirichlet boundary value in the half-line; see Proposition \ref{prop-nosol}. To complete the proof of
Theorem \ref{thm-l-d2} we will use the following Liouville theorem in dimension $n+1=2$.
\begin{prop} \label{chip} {\rm \bf (\cite{CCFS})}
Suppose that v weakly solves
\begin{equation}\label{eqn-chipotL}
\begin{cases}
-\Delta v\ge 0 &\text{in}\; \R^{2}_{+},\\
\frac{\partial v }{\partial \nu} \ge 0 & \text{on}\; \partial\R^{2}_{+}, \\
v \ge 0  &\text{in}\; \R^{2}_{+}.
\end{cases}
\end{equation}
Then, $v$ is a constant.
\end{prop}
As usual, very strong  Liouville theorems (but quite simple to
prove) hold in low dimensions, but not in higher ones. Compare (\ref{eqn-chipotL}) in low dimensions for supersolutions of the homogeneous linear problem
 with (\ref{eqn-li-ou}) for solutions of a precise nonlinear problem. The proof of Proposition \ref{chip}
 in \cite{CCFS} compared in an appropriate way the solution $v$ with
$\log (|\cdot|)$. For completeness, we give here  an alternative proof.

\medskip

\noindent{\bf Proof of Proposition \ref{chip}.} Replacing $v$ by $v-\inf_{\R^{2}_{+}}v\ge 0$, we may assume
$\inf_{\R^{2}_{+}}v=0$. Letting $w=1-v$, we have
\begin{equation}\label{eqn-w-one-v}
\begin{cases}
-\Delta w\le 0 &\text{in}\; \R^{2}_{+},\\
\frac{\partial w }{\partial \nu} \le 0 & \text{on}\; \partial\R^{2}_{+}, \\
w \le 1  &\text{in}\; \R^{2}_{+}.
\end{cases}
\end{equation}
In addition, $\sup_{\R^{2}_{+}}w=1$. Let $\xi_{R}\in C^{\infty}(\R^{2})$ be a function with compact
support in $B_{2R}(0)$, equal to $1$ in $B_{R}(0)$, and with
$|\nabla \xi_{R}|\le \frac{C}{R}$. Let
\[D_{R,2R}^{+}:=\{(x,y)\in\R^{2}\mid R\le |(x,y)|\le 2R, y
> 0\}.\]
Multiplying the first equation in (\ref{eqn-w-one-v}) by $w^{+}\xi_{R}^{2}$,
integrating in $\R^{2}_{+}$ and using the Neumann condition and $w^{+}\le 1$, we see that
\begin{align}\label{eqn-lv-pf01}
 \int_{\R^{2}_{+}}\xi_{R}^{2}|\nabla w^{+}|^{2}\le &
2\int_{D_{R,2R}^{+}}\xi_{R}\nabla \xi_{R}w^{+}\nabla w^{+}
\\\nonumber
\le & C\left(\int_{D_{R,2R}^{+}}|\nabla \xi_{R}|^{2}\right)^{1/2}
\left(\int_{D_{R,2R}^{+}} \xi_{R}^{2}|\nabla w^{+}|^{2}\right)^{1/2}\\
\le & C \left(\int_{D_{R,2R}^{+}} \xi_{R}^{2}|\nabla w^{+}|^{2}\right)^{1/2}.
\label{eqn-lv-pf03}
\end{align}
This leads, letting $R\uparrow \infty$, to $\int_{\R^{2}_{+}}|\nabla w^{+}|^{2} <\infty$. As a consequence
of this, the integral in (\ref{eqn-lv-pf03}) tends to zero as $R\rightarrow \infty$. Thus, by (\ref{eqn-lv-pf01}),
  \[\int_{\R^{2}_{+}}|\nabla
w^{+}|^{2}=0.\]
Thus, $w^{+}$ is constant, and since $\sup_{\R^{2}_{+}}w=1$, we conclude $w\equiv 1$. \hfill{$\Box$}

\begin{prop} \label{prop-sol-d2} Let $n\ge 2$ and
\[\R^{n+1}_{++}=\{z=(x_{1},x_{2},\cdots,x_{n},y)\mid x_{n}>0, y >
0\}.\]
 Assume that $v$ is a classical solution of
\begin{equation}\label{eqn-depend2}
\begin{cases}
\Delta  v=0 &  \text{in}\; \R^{n+1}_{++}, \\
v=0  &\text{on}\; \{x_{n}= 0, y\ge 0\},\\
\frac{\partial v}{\partial \nu}=v^{p} & \text{on}\; \{x_{n}>0, y=0\}, \\
 v>0 &  \text{in}\; \R^{n+1}_{++},
\end{cases}
\end{equation}
where  $1\le p\le \frac{n+1}{n-1}$.  Then, $v$ depends only on
$x_{n}$ and $y$.
\end{prop}
\noindent{\bf Proof.} We shall follow the steps of  \cite{GS}. Let
$e_{n}=(0,\cdots,0,1,0)$ and $N=n+1$. Consider the conformal
transformation
\[
\bar{z}=T(z)=\frac{z+e_{n}}{|z+e_{n}|^{2}}
\]
and the Kelvin transformation $w$ of $v$
\[
w(\bar{z})=|z+e_{n}|^{N-2}v(z)=|\bar{z}|^{2-N}v(z).
\]

Denote
$B_{1/2}^{+}(\frac{e_{n}}{2}):=\{\bar{z}=(\bar{x},\bar{y})\mid|\bar{z}-\frac{1}{2}
e_{n}|< \frac{1}{2}, \bar{y}> 0\}$,
$S_{1/2}^{+}(\frac{e_{n}}{2}):=\partial
B_{1/2}^{+}(\frac{e_{n}}{2})\cap \{\bar{y}>0\}$,
$\Gamma_{0,1/2}:=\partial B_{1/2}^{+}(\frac{e_{n}}{2})\cap
\{\bar{y}=0\}$.

Note that, through $T$, $\R^{n+1}_{++}=\{x_n>0,y>0\}$ gets mapped into
the half-ball $B_{1/2}^{+}(\frac{e_{n}}{2})$, the boundary
$\{x_{n}>0, y=0\}$ becomes the ball $\Gamma_{0,1/2}$, $\{x_{n}=
0, y\ge 0\}$ goes to the half-sphere $S_{1/2}^{+}(\frac{e_{n}}{2})$,
and the infinity goes to $\bar{z}=0$.

We see that $w$ satisfies
\begin{equation*}
\begin{cases}
\Delta  w=0 &  \text{in}\; B_{1/2}^{+}(\frac{e_{n}}{2}), \\
w=0 & \text{on}\; S_{1/2}^{+}(\frac{e_{n}}{2}), \\
\frac{\partial w(\bar{z})}{\partial \nu}=|\bar{z}|^{p(N-2)-N}w^{p}(\bar{z})  &\text{on}\; \Gamma_{0,1/2},\\
 w>0 &  \text{in}\; B_{1/2}^{+}(\frac{e_{n}}{2}).
\end{cases}
\end{equation*}
Since $|\bar{z}|^{p(n-1)-(n+1)}$ is nonincreasing in the
$\bar{z}_{i}$ direction for all $i=1,\cdots,n-1$ (in fact, in any
direction orthogonal to the $\bar{z}_{n}$-axis), the moving planes
method used as in \cite{CCFS} gives that $w$ is symmetric about all
the $\bar{z}_{i}$-axis for $i=1,\cdots,n-1$. This leads to
$w=w(|\bar{z}'|,\bar{z}_{n}, \bar{y})$, where
$\bar{z}'=(\bar{z}_{1},\cdots,\bar{z}_{n-1})$ and hence
$v=v(|x'|,x_{n},y)$. Now, since we may perform the Kelvin's
transform with respect to any point $(-x_{0}',-1,0)$ ---and not only
with respect to $x_{0}'=0$ as before---  we conclude
that $v=v(x_{n},y)$ as claimed. \hfill{$\Box$}

\begin{prop} \label{prop-nosol}
Assume that  $f$ is a $C^{1,\alpha}$ function for some $\alpha\in (0,1)$, such that $f>0$ in $(0,\infty)$  and $f(0)=0$. Let
 $C$ be a positive constant. Then there is no bounded solution of problem
\begin{equation}\label{eqn-2d}
\begin{cases}
\Delta v=0 &  \text{in}\;  \R^{2}_{++}=\{x>0,y>0\}, \\
v=0  &\text{on}\; \{x= 0, y\ge 0\}, \\
\frac{\partial v}{\partial \nu}=f(v) & \text{on}\; \{x>0, y=0\}, \\
0<v\le C &  \text{in  }\;  \R^{2}_{++}.
\end{cases}
\end{equation}
\end{prop}
\noindent {\bf Proof.} We use some tools developed in \cite{CS05}.

Suppose by contradiction that  there is such solution $v$. First, we
claim that $v(x,0)\rightarrow 0$ as $x\rightarrow \infty$. Suppose
by contradiction that there exists a sequence
$a_{m}\rightarrow\infty$ $(m\rightarrow \infty)$ such that
$v(a_{m},0)\rightarrow \alpha>0$. Let $v_{m}(x,y):=v(x+a_{m},y)$. It
is clear that $v_{m}$ is a solution of (\ref{eqn-2d}) in
$U_{m}:=\{(x,y)\mid x>-a_{m}, y>0\}$. Moreover,
$v_{m}(0,0)=v(a_{m},0)\rightarrow \alpha$. Therefore there exists a
subsequence, still denoted by $ v_{m}$, such that $ v_{m}\rightarrow
v$ in $C^{2}_{\rm loc}(\overline{\R_{+}^{2}})$ as $m\rightarrow
\infty$, and  $v$ is a solution of
\begin{equation}\label{eqn-2d-s}
\begin{cases}
\partial_{xx}v+  \partial_{yy}v=0 &  \text{in}\; \R^{2}_{+}=\{(x,y)\mid y>0\}, \\
\frac{\partial v}{\partial \nu}=f(v)\ge 0 & \text{on}\; \{ y=0\}, \\
0\le v\le C&  \text{in  }\;\{y>0\}.
\end{cases}
\end{equation}
Notice that
\[
v(0,0)=\alpha>0.
\]
On the other hand, by Proposition  \ref{chip} we know that $v$ is
identically constant. This is impossible due to the nonlinear
Neumann condition,  since $f>0$ in $(0,\infty)$ and
$f(v(0,0))=f(\alpha)>0$. We conclude the claim, that is,
$v(x,0)\rightarrow 0$ as $x\rightarrow +\infty$.

Note that we can reflect the function $v$ with respect to
$\{x=0,y>0\}$, $\tilde{v}(x,y)=-v(-x,y)$ for $x<0$, and obtain a
bounded harmonic function $\tilde{v}$ in all $\R^{2}_{+}=\{y>0\}$,
since $v\equiv 0$ on $\{x=0, y>0\}$. Applying interior gradient
estimates to the bounded harmonic function $\tilde{v}$ in the ball
$B_{t}(x,t)\subset \R^{2}_{+}$, we obtain
\[ |\nabla v (x,t)|\le \frac{C\|v\|_{\infty}}{t}
\le \frac{C}{t},\quad \text{for all}\; t>\frac{1}{2}, x>0.\] On the
other hand,
by the results of \cite{CS05} applied to the solution $\tilde{v}$ in $\R^{2}_{+}$
 (or equivalently by the proof of Proposition  \ref{prop-non-reg} of this paper; note that $f(0)=0$),
we have that $|\nabla v|$ and $|D^{2}v|$ are bounded in
$\overline{\R^{2}_{++}}\cap \{0\le y\le 1\}$.  We conclude that
$|\nabla v|$ and $|D^{2}v|$ are bounded in $\overline{\R^{2}_{++}}$
and
\[ |\nabla v (x,t)|\le  \frac{C}{t+1},\quad \text{for all}\; t>0, x>0.\]
Using interior estimates for harmonic functions as before, but now with
the partial derivatives of $v$ instead of $v$, it follows that
\[ |D^{2} v (x,t)|\le  \frac{C}{t^{2}+1},\quad \text{for all}\; t>0, x>0.\]
Moreover, we have
\[
\Big|\frac{\partial}{\partial
x}\left\{\frac{|\partial_{x}v(x,t)|^{2}-|\partial_{y}v(x,t)|^{2}}{2}\right\}\Big|\le
\frac{C}{t^{3}+1}.
\]
 By these facts, we  see that the function
\[
\Phi(x):=\int_{0}^{+\infty}\frac{|\partial_{x}v(x,t)|^{2}-|\partial_{y}v(x,t)|^{2}}{2}\,dt
\]
is well defined and $\frac{d\Phi}{dx}$ is also.

Using the $\lim\limits_{t\rightarrow \infty}|\nabla v(x,t)|=0$, we
obtain, for $F(v)=\int_{0}^{v}f(s)\,ds$,
\begin{equation*}
\begin{split}
\frac{d}{dx}&[\Phi(x)+F(v(x,0))]\\
&=
\int_{0}^{+\infty}[\partial_{xx}v\partial_{x}v-\partial_{y}v\partial_{xy}v](x,t)\,dt
+[f(v)\partial_{x}v](x,0)\\
&=[\partial_{y}v\partial_{x}v +f(v)\partial_{x}v](x,0)=0,
\end{split}
\end{equation*}
thanks to the harmonicity of $v$ and the Neumann boundary condition.
 This leads to the Hamiltonian-type identity
\[
\Phi(\cdot)+F(v(\cdot,0))\;\text{is identically constant in}\;(0,+\infty).
\]

Furthermore, using that $\lim\limits_{x\rightarrow +\infty}v(x,0)=0$, and
that $\lim\limits_{x\rightarrow +\infty}v(x,y)=~0$ uniformly
in compact sets in $y$ (we can prove this by the same previous argument leading to
$\lim\limits_{x\rightarrow +\infty}v(x,0)=0$), together with the above
bounds for $|\nabla v(x,y)|$ for $y$ large, we deduce
\[
\lim_{x\rightarrow+\infty}\Phi(x)=0.
\]
From all these we obtain
\[
\Phi(x)+F(v(x,0))\equiv 0, \;\mbox{for}\; x>0.
\]
Since $v=0$ and thus $\partial_{y}v=0$ along the $y$-axis, we see by
the definition of $\Phi(0)$ that
\[
0=\Phi(0)+F(v(0,0))=\Phi(0)=\frac{1}{2}\int_{0}^{+\infty}|\partial_{x}v|^{2}(0,t)\,dt.
\]
This implies that $\partial_{x}v=0$ on $\{x=0,y>0\}$, which
contradicts Hopf's lemma. Thus, the contradiction means that there is
no positive bounded solution of the problem. \hfill{$\Box$}

\medskip

Before  proving Theorems \ref{thm-l-d2} and  \ref{teo-ap}, let us make some comments.
\begin{remark}{\rm
 Theorem \ref{thm-l-d2} is still open without the boundedness assumption on $v$.

In this respect, let us give some examples of problems in the
quarter plane $\R^{2}_{++}$. The function $v(x,y)=x$ is
an unbounded solution of problem
\begin{equation*}
\begin{cases}
-\Delta v=0,\; v\ge 0  &\mbox{in}\;\;\R^{2}_{++},\\
v=0 & \mbox{on}\; \{x=0, y>0\},\\
\frac{\partial v}{\partial \nu}= 0 & \mbox{on}\; \{x>0, y=0\}.
\end{cases}
\end{equation*}
This tells us that the result of  Proposition \ref{chip} (which
did not require boundedness of the solution in the half-plane) does not
hold in the quarter plane.

 On the other hand, it is clear that
 $v(x,y)=\frac{\pi}{2}\arctan \frac{x}{y+1}$  satisfies $\Delta v=0$ and
$-\partial_{y}v\mid_{y=0}=\frac{\pi x}{2(1+x^{2})}\ge 0$ for $x>0$. Hence,
there exists a bounded harmonic function in the quarter plane
$\R^{2}_{++}$ such that
\begin{equation*}
\begin{cases}
-\Delta v=0,\; v\ge 0  &\mbox{in}\;\;\R^{2}_{++},\\
v=0 & \mbox{on}\; \{x=0, y>0\},\\
\frac{\partial v}{\partial \nu}\ge 0 & \mbox{on}\; \{x>0, y=0\}.
\end{cases}
\end{equation*}
Thus the nonlinear condition $\frac{\partial v}{\partial \nu}=v^{p}$
on $\{y=0\}$ is  important in Theorem \ref{thm-l-d2}.}
 \end{remark}

\noindent{\bf Proof of Theorem \ref{thm-l-d2}.} It follows from
Propositions \ref{prop-sol-d2} and \ref{prop-nosol}.
\hfill{$\Box$}

\medskip

\noindent{\bf Proof of Theorem \ref{teo-ap}.} We know by Theorem
\ref{prop-teo-bk} and Proposition \ref{prop-non-reg}  that  all weak
solutions $u $ of (\ref{eqn-frac}) belong to $
C^{2}(\overline{\Omega})\cap C_{0}(\overline{\Omega})$.
 Assume by contradiction that the theorem is not true and hence
 that there is a sequence $u_{m}$ of solutions  of (\ref{eqn-frac})
 with
\[
K_{m}=\|u_{m}\|_{L^{\infty}(\Omega)}\rightarrow \infty.
\]

 Since $v_{m}=\text{h-ext}(u_{m})$ is a positive harmonic function in $\C$ vanishing on
$\partial_{L}\C$, we have that $v_{m}$ has also $K_{m}$ as maximum
in $\C$ and that it is attained  at a point $(x_{m},0)\in
\Omega\times\{0\}$. Let
\[\Omega_{m}=K_{m}^{p-1}(\Omega-x_{m})\]
 and define
\[
\tilde{v}_{m}(x,y)=K_{m}^{-1}v(x_{m}+K_{m}^{1-p}x,K_{m}^{1-p}y),\quad
x\in \Omega_{m}, y>0.
\]
We have that
$\|\tilde{v}_{m}\|_{L^{\infty}(\Omega_{m}\times(0,\infty))}\le 1$
and
\begin{equation}\label{eqn-secap-K}
\begin{cases}
\Delta\tilde{v}_{m}=0,\; & \mbox{in}\;\;\C_{m}:=\Omega_{m}\times (0,\infty),\\
\tilde{v}_{m}=0 &\mbox{on}\;\;\partial
\Omega_{m}\times(0,\infty),\\
\frac{\partial\tilde{v}_{m}}{\partial \nu}=\tilde{v}_{m}^{p}&\mbox{on}\;\Omega_{m}\times\{0\},\\
\tilde{v}_{m}>0& \mbox{in}\;\;\C_{m}.
\end{cases}
\end{equation}
Notice that \[\tilde{v}_{m}(0,0)=1.\]

Let
 \[
 d_{m}=\text{dist}(x_{m},\partial \Omega).
\]
Two cases may occur as $m\rightarrow \infty$; either case (a):
\[
K_{m}^{p-1}d_{m}\rightarrow \infty
\]
for a subsequence still denoted as before, or case (b):
\[
K_{m}^{p-1}d_{m}\;\text{ is bounded}.
\]

 If case (a) occurs, we have that $B_{K_{m}^{p-1}d_{m}}(0)
 =K_{m}^{p-1}B_{d_{m}}(0)\subset \Omega_{m}$
 and that $K_{m}^{p-1}d_{m}\rightarrow \infty$.
 By local compactness (Arzel\`{a}-Ascoli) of bounded solutions to (\ref{eqn-secap-K})
 (recall $\|\tilde{v}_{m}\|_{L^{\infty}(\Omega_{m})}
 \le 1$), through a subsequence, we obtain a solution $\tilde{v}$ of problem (\ref{eqn-li-ou})
 in all of $\R^{n+1}_{+}=\R^{n}\times(0,\infty)$
---note that $\tilde{v}_{m}(0,0)=1$ leads to  $\tilde{v}(0,0)=1$ and hence $\tilde{v}\not\equiv 0$ \& $\tilde{v}>0$. This is a contradiction to Theorem \ref{li-ou}.

Assume now that case (b), $K_{m}^{p-1}d_{m}$ is bounded, occurs.
Note first that since the right-hand side of problem (\ref{eqn-bd1}) for $v_{m}$ satisfies $|v_{m}|^{p}=v_{m}^{p}\le
K_{m}^{p}$, we deduce from the proofs of Proposition \ref{prop-linear-reg} (iii) and (iv) that $\|\nabla u_{m}\|_{L^{\infty}(\Omega)}\le
CK^{p}_{m}$ for a constant $C$ independent of $m$. Now, since
$u_{m}\!\!\mid_{\partial\Omega}\equiv 0$ (where $u_{m}=v_{m}(\cdot,0)$), we get
\[K_{m}=v_{m}(x_{m},0)\le \|\nabla u_{m}
\|_{L^{\infty}(\Omega)}\text{dist}(x_{m},\partial \Omega)\le
CK_{m}^{p}d_{m}.\]
 We deduce that
\[
0<c\le K_{m}^{p-1}d_{m}
\]
for some positive constant $c$. Thus, in this case (b), we may
assume that, up to a subsequence,
\begin{equation}\label{eqn-pf-ap-last}
K_{m}^{p-1}d_{m}\rightarrow a\in (0,\infty)
\end{equation}
for some constant $a>0$.

We deduce that, up to a certain rotation of $\R^{n}$ for each index
$m$, since we have (\ref{eqn-pf-ap-last}), $K_{m}^{p-1}\rightarrow \infty$, $d_{m}\rightarrow
0$, and that $B_{K_{m}^{p-1}d_{m}}(0)$ is tangent to $\partial
\Omega_{m}$, the domains $\Omega_{m}$ converge to the half-space
$\R^{n}_{+}=\{x_{n}>-a\}$. Thus, through a subsequence of
$\tilde{v}_{m}$, we obtain a solution $\tilde{v}$ of problem
(\ref{eqn-depend2}) in $\R^{n+1}_{++}=\{x_{n}>-a,y>0\}$ with
$\tilde{v}$ bounded by $1$ and $\tilde{v}>0$ (since
$\tilde{v}_{m}(0,0)=1$ for all $m$).
 This is a contradiction with Theorem~\ref{thm-l-d2}.  \hfill{$\Box$}

\begin{remark}{\rm From Theorem  \ref{teo-ap} we have
 a priori bounds for solutions of  problem (\ref{eqn-frac})
  with $f(u)=u^{p}$, $1<p<\frac{n+1}{n-1}$.
As a consequence, by using blow-up techniques and topological degree theory, 
one can obtain existence
of positive solutions for related problems ---for instance, for nonlinearities $f(x,u)$
of power type, as well as other boundary conditions. See Gidas-Spruck~\cite{GS}
for some of these applications when the operator is the classical Laplacian.
}
\end{remark}

\setcounter{equation}{0}
\section{Symmetry of solutions}\label{sec-h-sys}
The goal of this section is to prove a symmetry result of
Gidas-Ni-Nirenberg type for positive solutions of nonlinear problems
involving the operator $A_{1/2}$, as stated in Theorem
\ref{teo-sys}, by using the moving planes method. For this, we work with the equivalent local problem
 (\ref{eqn-frac2}) and derive the following.

 \begin{teo}\label{teo-v-sys7}
Assume that $\Omega$ is a bounded smooth domain of  $\R^{n}$ which
is convex in the $x_{1}$ direction and symmetric with respect to the
hyperplane $\{x_{1}=0\}$. let  $f$ be Lipschitz continuous and let  $v\in C^{2}(\overline{\C})$ be a solution of
{\rm(\ref{eqn-frac2})}, where $\C=\Omega\times (0,+\infty)$. Then, $v$ is symmetric with respect to
$x_{1}$, i.e., $v(-x_{1},x',y)=v(x_{1},x',y)$ for all
$(-x_{1},x',y)\in \C$. In addition, $\frac{\partial v}{\partial
x_{1}}<0$ for $x_{1}>0$.
\end{teo}

 \noindent {\bf Proof of Theorems \ref{teo-sys} and \ref{teo-v-sys7}.}
 It suffices to prove Theorem \ref{teo-v-sys7}. From it, Theorem \ref{teo-sys} follows
 immediately.

 Let
$x=(x_{1},x')\in \Omega$ and $\lambda>0$. Consider the sets
\[
\Sigma_{\lambda}=\{(x_{1},x')\in \Omega\mid x_{1}>\lambda\}
\quad\mbox{and}\quad T_{\lambda}=\{(x_{1},x')\in \Omega\mid
x_{1}=\lambda\}.
\]
For $x\in \Sigma_{\lambda}$, define
$x_{\lambda}=(2\lambda-x_{1},x')$. By the hypotheses on the domain~$\Omega$ we see that
\[
\{x_{\lambda}\mid x\in \Sigma_{\lambda}\}\subset \Omega.
\]

Recall that $v\in C^{2}(\overline{\C})$  is a solution of
\begin{equation*}
\left\{
\begin{array}{ll}
 \Delta v=0 &  \mbox{in}\quad {\mathcal C}=\Omega\times (0,\infty),\\
 v=0 & \mbox{on}\quad \partial_{L}{\mathcal C}=\partial\Omega\times [0,\infty),\\
\frac{\partial v}{\partial \nu}=f(v)& \mbox{on} \quad \Omega\times \{0\},\\
v>0 &  \mbox{in}\quad {\mathcal C}.
\end{array}
\right.
\end{equation*}
For $(x,y)\in\Sigma_{\lambda}\times[0,\infty)$, let us define
\[
v_{\lambda}(x,y)=v(x_{\lambda},y) =v(2\lambda-x_{1},x',y)
\]
 and
\[w_{\lambda}(x,y)=(v_{\lambda}-v)(x,y).\]
Note that $v_{\lambda}$ satisfies
\begin{equation*}
\left\{
\begin{array}{ll}
 \Delta v_{\lambda}=0 &  \mbox{in}\quad \Sigma_{\lambda}\times (0,\infty),\\
 v_{\lambda}\ge 0 & \mbox{on}\quad (\partial\Omega \cap \overline{\Sigma}_{\lambda})\times (0,\infty),\\
\frac{\partial v_{\lambda}}{\partial \nu}=f(v_{\lambda})& \mbox{on}
\quad \Sigma_{\lambda}\times \{0\}.
\end{array}
\right.
\end{equation*}
Thus, since $\partial \Sigma_{\lambda}= (\partial
\Omega\cap\overline{\Sigma}_{\lambda})\cup T_{\lambda}$ and
$w_{\lambda}\equiv 0$ on $T_{\lambda}$, we have that $w_{\lambda}$
satisfies
\begin{equation}\label{eqn-pf-sys-cp}
\left\{
\begin{array}{ll}
 \Delta w_{\lambda}=0 &  \mbox{in}\quad \Sigma_{\lambda}\times (0,\infty),\\
 w_{\lambda}\ge 0 & \mbox{on}\quad (\partial \Sigma_{\lambda})\times (0,\infty),\\
\frac{\partial w_{\lambda}}{\partial
\nu}+c_{\lambda}(x)w_{\lambda}=0& \mbox{on} \quad
\Sigma_{\lambda}\times \{0\},
\end{array}
\right.
\end{equation}
where \[c_{\lambda}(x,0)=-\frac{f(v_{\lambda})-f(v)}{v_{\lambda}-v}.\]
Note that $c_{\lambda}(x,0)\in L^{\infty}(\Sigma_{\lambda})$.

Let $\lambda^{*}=\sup \{ \lambda\mid \Sigma_{\lambda}\neq
\emptyset\}$ and let $\varepsilon>0$ be a small number. If
$\lambda\in (\lambda^{*}-\varepsilon,\lambda^{*})$, then
$\Sigma_{\lambda}$ has small measure and we have, by part (ii) of
Proposition~\ref{lem-max3} (applied with $\Omega$ replied  with
$\Sigma_{\lambda}$), that
\[
w_{\lambda}\ge 0\quad \mbox{in}\;\Sigma_{\lambda}\times(0,\infty).
\]
Note here that $\Sigma_{\lambda}$ is not a smooth domain but that part (ii) of
Proposition \ref{lem-max3}  does not require smoothness of the domain. 
By the strong maximum principle, Lemma \ref{lem-max2}, for problem
(\ref{eqn-pf-sys-cp}) we see that $w_{\lambda}$ is
identically equal to zero or strictly positive in
$\Sigma_{\lambda}\times(0,\infty)$. Since $\lambda>0$, we have
$w_{\lambda}
>0$ in $(\partial \Omega\cap \partial \Sigma_{\lambda})\times (0,\infty)$,
and hence we conclude that
$w_{\lambda}
>0$ in $\Sigma_{\lambda}\times (0,\infty)$.

Let $\lambda_{0}=\inf \{ \lambda>0\mid
 w_{\lambda}\ge 0 \;\text{in}\; \Sigma_{\lambda}\times (0,\infty)\}$. We are going to prove that
$\lambda_{0}=0$. Suppose that $\lambda_{0}>0$ by contradiction.
First, by continuity, we have $w_{\lambda_{0}}\ge 0$ in
$\Sigma_{\lambda_{0}}\times(0,\infty)$. Then, as before, we deduce
$w_{\lambda_{0}}> 0$ in $\Sigma_{\lambda_{0}}\times (0,\infty)$.
Next, let $\delta>0$ be a constant and $K\subset
\Sigma_{\lambda_{0}}$ be a compact set such that
$|\Sigma_{\lambda_{0}}\setminus K|\le \delta/2$. We have
$w_{\lambda_{0}}(\cdot,0)\ge \eta>0$ in $K$ for some constant
$\eta$, since $K$ is compact. Thus, we obtain that
$w_{\lambda_{0}-\varepsilon}(\cdot,0)>0$ in $K$ and that
$|\Sigma_{\lambda_{0}-\varepsilon}\setminus K|\le \delta$ for
$\varepsilon$ small enough.

Now we apply again part (ii) of Proposition~\ref{lem-max3}
in $\Sigma_{\lambda_{0}-\varepsilon}\times(0,\infty)$  to the function $w_{\lambda_{0}-\varepsilon}$.
 We know that $w_{\lambda_{0}-\varepsilon}(\cdot,0)\ge 0$ in $K$, and hence
 $\{w_{\lambda_{0}-\varepsilon}<0\}\subset
\Sigma_{\lambda_{0}-\varepsilon}\setminus K$, which has measure at
most $\delta$. We take~$\delta$ to be the constant of part (ii) of
Proposition~\ref{lem-max3}. We deduce that
\[
w_{\lambda_{0}-\varepsilon}\ge
0\quad\mbox{in}\;\Sigma_{\lambda_{0}-\varepsilon}\times(0,\infty).
\]
This is a contradiction to the definition of $\lambda_{0}$. Thus,
$\lambda_{0}=0$.

We have proved, letting $\lambda\downarrow \lambda_{0}=0$ that
\[
v(-x_{1},x',y)\ge v(x_{1},x',y) \quad \mbox{in}\;
(\Omega\cap\{x_{1}>0\})\times(0,\infty)
\]
and, since $w_{\lambda}=0$ on $T_{\lambda}$,
\[
\partial_{x_{1}}v=-\frac{1}{2}\frac{\partial w_{\lambda}}{\partial
x_{1}}< 0\quad\mbox{for}\;x_{1}>0,
\]
by Hopf's lemma. Finally replacing $x_{1}$ by $-x_{1}$, we deduce the
desired symmetry $v(-x_{1},x',y)=v(x_{1},x',y)$.
 \hfill{$\Box$}

\bigskip

\noindent {\bf Acknowledgements:} Both authors were supported by Spain Government
grants MTM2005-07660-C02-01 and MTM2008-06349-C03-01. The second author was supported by CONICYT Becas
de Postgrado of Chile and the Programa de Recerca del Centre de Recerca Matem\`{a}tica, Barcelona, Spain.

\end{document}